\documentclass[final,3p,times]{elsarticle}
\usepackage{mathptmx}
\usepackage{amsfonts}
\usepackage{amsmath}
\usepackage{amssymb}
\usepackage[most]{tcolorbox}

\usepackage{hyperref}
\usepackage{xparse}
\usepackage{hhline}
\usepackage[author={Yous van Halder}]{pdfcomment}
\def\exampletext{Example} 

\NewDocumentEnvironment{testexample}{ O{} }
{
\colorlet{colexam}{black} 
\newtcolorbox[use counter=testexample]{testexamplebox}{%
    empty,
    title={\exampletext: #1},
    attach boxed title to top left,
       minipage boxed title,
    boxed title style={empty,size=minimal,toprule=0pt,top=4pt,left=3mm,overlay={}},
    coltitle=colexam,fonttitle=\bfseries,
    before=\par\medskip\noindent,parbox=false,boxsep=0pt,left=3mm,right=0mm,top=2pt,breakable,pad at break=0mm,
       before upper=\csname @totalleftmargin\endcsname0pt, 
    overlay unbroken={\draw[colexam,line width=.5pt] ([xshift=-0pt]title.north west) -- ([xshift=-0pt]frame.south west); },
    overlay first={\draw[colexam,line width=.5pt] ([xshift=-0pt]title.north west) -- ([xshift=-0pt]frame.south west); },
    overlay middle={\draw[colexam,line width=.5pt] ([xshift=-0pt]frame.north west) -- ([xshift=-0pt]frame.south west); },
    overlay last={\draw[colexam,line width=.5pt] ([xshift=-0pt]frame.north west) -- ([xshift=-0pt]frame.south west); },%
    }
\begin{testexamplebox}}
{\end{testexamplebox}\endlist}

\providecommand{\mbf}[1]{\mathbf{#1}}

\providecommand{\norm}[1]{\left\Vert #1\right\Vert}
\providecommand{\der}[2]{\frac{\partial #1}{\partial #2}}

\providecommand{\vecbrace}[1]{\left(\begin{array}{c} #1 \end{array}\right)}

\newdefinition{rmk}{Remark}
\newproof{pf}{Proof}
\newproof{pot}{Proof of Theorem \ref{thm2}}

\journal{Journal of Computational Physics}

\begin{document}
\begin{frontmatter}



\title{An adaptive minimum spanning tree multi-element method for \\
uncertainty quantification of smooth and discontinuous responses}


\author[main]{Y. van Halder}
\address[main]{Centrum Wiskunde \& Informatica (CWI), Science Park 123, 1098 XG, Amsterdam, the Netherlands}
\ead{y.van.halder@cwi.nl}
\author[main]{B. Sanderse}
\cortext[cor1]{Corresponding author}
\author[barry]{B. Koren}
\address[barry]{Eindhoven University of Technology, P.O. Box 513, 5600 MB, Eindhoven, the Netherlands}

\begin{abstract}
A novel approach for non-intrusive uncertainty propagation is proposed. Our approach overcomes the limitation of many traditional methods, such as generalised polynomial chaos methods, which may lack sufficient accuracy when the quantity of interest depends discontinuously on the input parameters. As a remedy we propose an adaptive sampling algorithm based on minimum spanning trees combined with a domain decomposition method based on support vector machines. The minimum spanning tree determines new sample locations based on both the probability density of the input parameters and the gradient in the quantity of interest. The support vector machine efficiently decomposes the random space in multiple elements, avoiding the appearance of Gibbs phenomena near discontinuities. On each element, local approximations are constructed by means of least orthogonal interpolation, in order to produce stable interpolation on the unstructured sample set. The resulting minimum spanning tree multi-element method does not require initial knowledge of the behaviour of the quantity of interest and automatically detects whether discontinuities are present. We present several numerical examples that demonstrate accuracy, efficiency and generality of the method.
\end{abstract}

\end{frontmatter}
{

\section{Introduction}
\noindent Uncertainty quantification (UQ) has become increasingly important for complex engineering applications. Determining and quantifying the influence of parametric and  model-form uncertainties is essential for a wide range of applications: from turbulent flow phenomena \cite{xiao_quantifying_2016, edeling_simplex-stochastic_2016}, aerodynamics \cite{simon_gpc-based_2010}, biology \cite{cho_experimental_2003, abagyan_biased_1994} to design optimisation \cite{constantinescu_computational_2011, mateos_monte_2000, papadrakakis_reliability-based_2002}. We are interested among others in liquid impact problems \cite{hopfinger_liquid_nodate,marsooli_3-d_2014,larocque_3-d_2013}.

For problems which have a complex underlying model, one often uses so-called \textit{non-intrusive} methods, also known as \textit{sampling} methods. The model is solved deterministically a number of times, and a stochastic solution is constructed using these deterministic samples. A well known sampling method for propagating uncertainties through a model is the \textit{Monte Carlo} method \cite{hammersley_monte_2013}. Despite its easy implementation and wide applicability, the Monte Carlo method suffers from slow convergence with increasing number of model evaluations, when approximating the Quantity of Interest (QoI). As a consequence of this slow convergence rate, many samples are required for obtaining high quality stochastic solutions. Therefore, Monte Carlo methods are not suitable for problems with a complex underlying model. As an alternative to Monte Carlo methods, \textit{Stochastic Collocation} (SC) methods \cite{xiu_fast_2009, babuska_stochastic_2007, nobile_sparse_2008} were introduced, replacing the slow convergence of Monte Carlo by an exponential convergence rate. The introduction of SC methods resulted in a decrease of required samples to achieve a certain accuracy in comparison to Monte Carlo methods. For a smooth QoI as a function of the uncertainties, fast convergence is achieved. However, if the QoI is highly non-linear or discontinuous, \textit{Gibbs phenomena} \cite{canuto_spectral_2012} may occur, which deteriorate the accuracy globally. To avoid the occurrence of Gibbs phenomena, several alternatives to the SC methods were introduced \cite{witteveen_adaptive_2009, jakeman_local_2011}, but they focus solely on discontinuous QoIs, leading to a significant increase in the number of samples needed for approximating smooth QoIs. One interesting method is the \textit{Multi-Element Stochastic Collocation} method (ME-SC) \cite{wan_multi-element_2006, foo_multi-element_2008}. The idea of ME-SC is to decompose the domain, spanned by the uncertainties, into smaller non-overlapping elements, in each of which the QoI is amenable for using an SC method. Gibbs phenomena still appear in the elements where there is a discontinuity in the QoI, but they are confined to these specific elements. Improving the multi-element approach is an active field of research and focuses on more efficient and robust domain decomposition. Jakeman et al. \cite{jakeman_minimal_2013} proposed the minimal multi-element method, which uses discontinuity detection based on polynomial annihilation to detect discontinuities, and divides the domain along the discontinuities. As a result, the Gibbs phenomena are removed completely. Although this discontinuity detection algorithm is accurate, the total number of samples needed to determine the discontinuity location is still too high if sampling the model is expensive. Therefore, a discontinuity detection algorithm which performs well for a lower number of samples was proposed by Gorodetsky et al. \cite{gorodetsky_efficient_2014}. This discontinuity detection uses polynomial annihilation in combination with \textit{Support Vector Machines} (SVMs) to divide the domain into elements along the so-called SVM classification boundary. Even though both discontinuity detection algorithms \cite{jakeman_minimal_2013,gorodetsky_efficient_2014} perform well, using them for approximating smooth QoIs can become prohibitively expensive, as both methods focus solely on finding the discontinuities. It is often unknown in advance if a QoI is smooth or discontinuous and choosing a method which is suited for one of both only is not recommended. Our goal is to create a surrogate model that works for both smooth and discontinuous QoIs and which requires no initial knowledge about the QoI.

A novel domain decomposition method in combination with adaptive sampling of the QoI is used for constructing this surrogate. The adaptive sampling procedure in our method is based on \textit{Minimum Spanning Trees} (MST) \cite{kruskal_shortest_1956, graham_history_1985}, which add new sample points at places which are associated with a high probability density and/or where the QoI changes rapidly. The adaptively placed samples are classified and an SVM \cite{gorodetsky_efficient_2014, el-tawel_edge_2015, varma_pixel-based_2016, burges_tutorial_1998, scholkopf_learning_2001} is used to obtain a classification boundary, which serves as an approximation for the discontinuity location. The decomposition of the random space in this way leads to elements on which each local QoI is amenable for interpolation without Gibbs phenomena. For constructing a surrogate model in each element a \textit{least orthogonal interpolant} \cite{narayan_stochastic_2012} is employed, which is suited for interpolation on the scattered data points that we obtain with our adaptive sampling. Our proposed method is abbreviated as MST-ME (Minimum Spanning Tree Multi-Element) method and is designed mainly for the purpose of uncertainty propagation. However, when assuming uniformly distributed probability density functions for the input parameters, the MST-ME method is also well suited to obtain a parametric solution of the partial differential equation (PDE) under consideration. This will also be illustrated in this paper.

This paper is outlined as follows: section \ref{sec:PD} briefly introduces the problem. Section \ref{sec:Algorithm} introduces the MST-ME method in detail. Finally, section \ref{sec:Results} demonstrates efficiency and accuracy of our method when applied to analytical test-cases. Complex test-cases related to sloshing impact problems, i.e., shallow water dam break and 3D dam break, are also discussed.
\section{Problem Description}
\label{sec:PD}
\noindent Quantifying the effects of uncertainties in computational engineering typically consists of three steps: (i) the input uncertainties are characterised in terms of a probability density function (PDF), which follows from observations or physical evidence; (ii) the uncertainties are propagated through the model; (iii) the outputs are post-processed, where the QoI is expressed in terms of its statistical properties. In the present work we focus on the propagation step, and the input distributions are assumed to be given. The goal is to solve the following stochastic problem:
\begin{equation}
\mathcal{L}(u ; \mbf{z}) = 0\ ,
\label{eq1:ModelProblem}
\end{equation}
where $\mbf{z} = (z_1, ..., z_d)\in I_{\mbf{z}}$ is the $d$-dimensional vector containing uncertain inputs, $u = u(\mbf{z})$ the QoI and $\mathcal{L}$ an operator, which represents the model. The operator $\mathcal{L}$ can be a (non-)linear partial differential operator, or any mathematical model that relates input $\mbf{z}$ to the QoI $u(\mbf{z})$. Both continuous QoIs $u\in C^0(I_{\mbf{z}})$, and discontinuous QoIs $u\not\in C^0(I_{\mbf{z}})$, are considered in this paper. The support set $I_{\mbf{z}}$ of the uncertain inputs $\mbf{z}$ is referred to as the \textit{random space}. The uncertain inputs are assumed to be characterised by a joint PDF $\rho(\mbf{z})$. The stochastic problem is solved non-intrusively by sampling the model \eqref{eq1:ModelProblem} at different locations $\mbf{z}_i$ in the random space, i.e.:
\begin{equation}
\mathcal{L}(u_i ; \mbf{z_i}) = 0\ ,
\end{equation}
where $u_i$ is the sampled solution at the collocation node $\mbf{z}_i$. Since the sampling is non-intrusive, black-box solvers for the operator $\mathcal{L}$ can be used. In this paper we are interested in finding the entire functional relation $u$ as a function of the uncertainties $\mbf{z}$, in terms of a surrogate model. A surrogate model $\tilde{u}$ of $u$ is constructed by interpolation, such that:
\begin{equation}
\tilde{u}(\mbf{z})\approx u(\mbf{z}),\ \ \text{for all}\ \ \mbf{z}\in I_{\mbf{z}}\ .
\end{equation}
When $u$ is smooth, it is possible to construct an approximation $\tilde{u}$ which converges exponentially fast to the exact solution $u$. However, if the QoI exhibits highly non-linear or discontinuous behaviour, then the accuracy of the approximation deteriorates globally, due to Gibbs phenomena. Multi-element methods divide the random space $I_{\mbf{z}}$ into a set of $N_E$ smaller elements $E_i$, such that the negative impact of the Gibbs phenomena is confined to a limited number of elements surrounding the discontinuity. The elements $E_i$ are non overlapping and span the entire random space, i.e.:
\begin{equation}
\cup_{i=1}^{N_E} E_i = I_{\mbf{z}}\ \ \text{and}\ \ E_i\cap E_j = \emptyset\ ,\ \ \text{if}\ \ i\neq j\ . 
\label{eq1:Elements}
\end{equation}
A local surrogate $\tilde{u}^{(i)}$ is constructed in each $E_i$:
\begin{equation}
\tilde{u}^{(i)}(\mbf{z})\approx u(\mbf{z}),\ \ \text{for all}\ \ \mbf{z}\in E_i\ .
\end{equation}
The global surrogate is given by patching the local surrogate models:
\begin{equation}
\tilde{u}(\mbf{z}) = \sum_{i=1}^{N_E} \tilde{u}^{(i)}(\mbf{z}) \mathcal{I}_{E_i}(\mbf{z})\ ,
\label{eq1:globalsolution}
\end{equation}
where $\mathcal{I}_{E_i}(\mbf{z})$ is the indicator function satisfying $\mathcal{I}_{E_i}(\mbf{z})=1$, if $\mbf{z}\in E_i$ and 0 otherwise. Standard multi-element methods utilise a tensor construction of hypercubes for defining the elements $E_i$ \cite{wan_multi-element_2006, foo_multi-element_2008}. While such a tensor construction removes the global effect of Gibbs phenomena, they can still appear locally in several elements.

An example of the standard tensorised multi-element approach for the approximation of a 2D function is shown in figure \ref{fig:MEExample}. Gibbs phenomena appear in the elements $E_{3,4}$, where a discontinuity is present in the exact function.
\begin{figure}[!h]
\centering
\includegraphics[width = 0.8\textwidth]{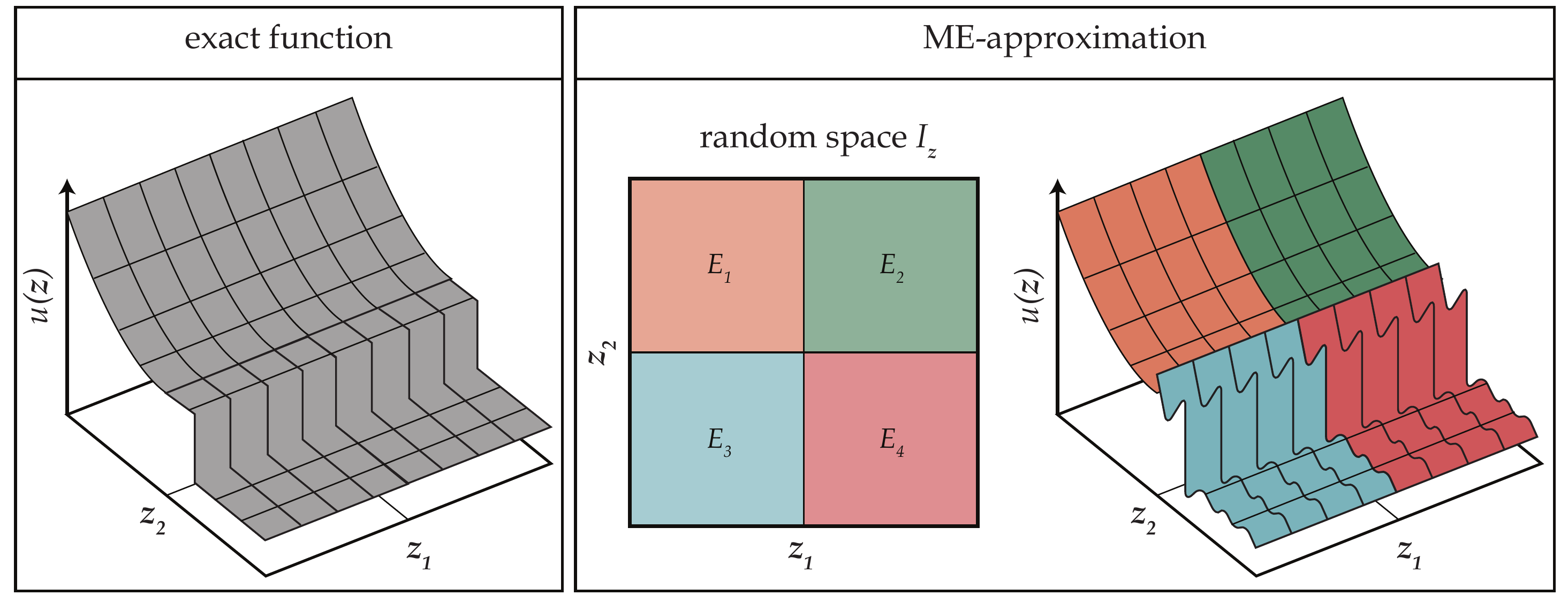}
\caption{\label{fig:MEExample} Example of the standard multi-element approach for approximating a 2D function exhibiting a discontinuity.}
\end{figure}

{\renewcommand{\thesubsection}{\Roman{subsection}}
\section{Minimum Spanning Tree Multi-Element Method}
\label{sec:Algorithm}
\noindent The art of an efficient multi-element stochastic collocation method lies in the choice of sampling points $\mbf{z}_i$, the choice of the elements $E_i$, and the reconstruction of the local approximations $\tilde{u}^{(i)}$. These are the focuses of this paper. Accordingly, the MST-ME method introduced here, is divided into three stages: (I) sampling, (II) domain decomposition and (III) local approximation construction:\\

\noindent{\indent I (\textit{choice of $z_i$}) Adaptive sampling of the QoI while taking into account both smooth and discontinuous regions, and the underlying PDF.\\
\indent II (\textit{choice of $E_i$}) Division of the random space into a minimal number of elements, such that the QoI is smooth within each element.\\
\indent III (\textit{construction of $\tilde{u}^{(i)}$}) Interpolation of the samples, while producing a stable interpolant.}\\

\noindent The QoI is adaptively sampled (I) by taking into account both the PDF and the QoI gradient information. The samples are distributed among different classes, such that within each class the local QoI is smooth. The classified samples are given as input to the domain classification algorithm (II). Instead of using a tensor based domain decomposition, the random space is divided into a minimal number of elements, in which the QoI is amenable for interpolation without Gibbs phenomena. This domain decomposition methodology was already introduced in \cite{jakeman_minimal_2013}, but in that work the number of samples required for determining proper elements was too high. In contrast, our domain decomposition method (II) uses the samples from the sampling algorithm (I) for determining proper elements, without the need to perform additional sampling. Local approximations (III) are constructed in each element, by using the least orthogonal interpolation method \cite{narayan_stochastic_2012}. The global approximation, the surrogate model, is given by the patched local approximations \eqref{eq1:globalsolution}.

Our method distinguishes itself from other methods, such as \cite{jakeman_local_2011, jakeman_minimal_2013}, by detecting if a QoI is smooth or discontinuous automatically, while not emphasising solely on either the discontinuity or the smooth regions. The nodes are placed by considering both smooth and discontinuous characteristics of the QoI, which results in sample locations that resolve both smooth and discontinuous regions. The combination of these sample locations and an SVM leads to accurate discontinuity detection, without placing samples in the random space where they are not needed.
\subsection{Sampling Algorithm}
\label{sec:SamplingProcedure}
\noindent The main idea behind our adaptive sampling algorithm is that we want to refine our surrogate model based on the QoI behaviour and the associated PDF of the random input variables. We achieve this by creating a graph that links the samples, and by assigning weights to edges of this graph, constructing an MST, and then adding samples on the most important edges of this MST.

The methodology is explained for a 2D QoI $u(z_1, z_2)$, but can be generalised easily to a higher dimensional QoI.\\

\noindent\textit{{Initial sample placement}}\\
The sampling procedure is started by placing initial sample points. Straightforward choices of the initial sample locations are shown in figure \ref{fig:InitialSamples}. These initial sample configurations are extensible to high-dimensional random spaces. Initial sample placement in both figure \ref{fig:InitialSamples}a and figure \ref{fig:InitialSamples}b introduces an anisotropy in the placement of subsequent samples. The initial sample placement in figure \ref{fig:InitialSamples}c is a good trade-off between the number of initial samples and the isotropy of subsequent samples.
\begin{figure}[!h]
\centering
\includegraphics[width = \textwidth]{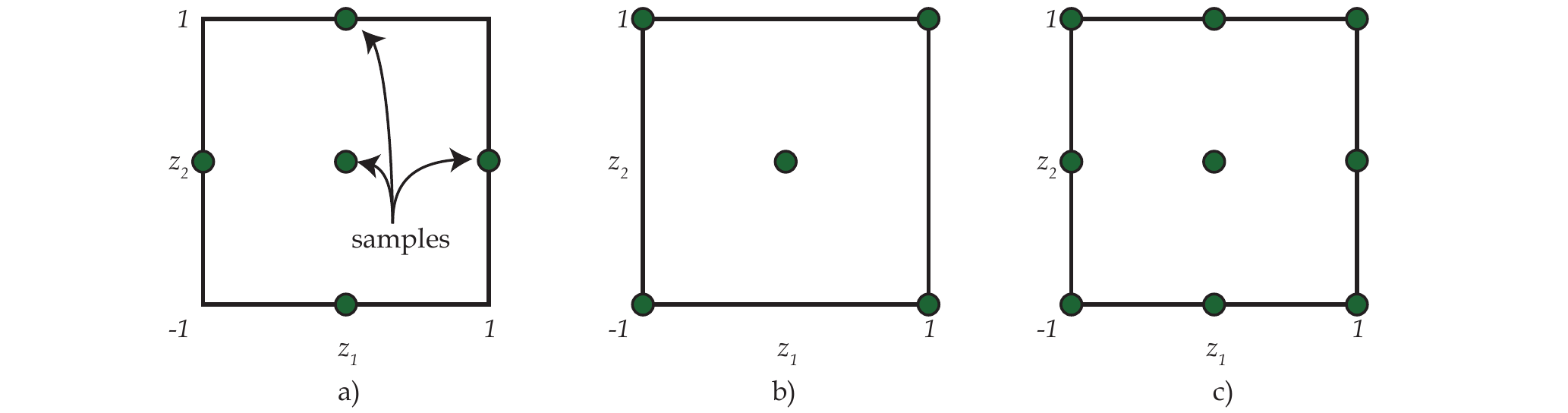}
\caption{\label{fig:InitialSamples} Initial sample location configurations. a) One sample in the middle and one sample at the face centres. b) One sample in the middle and one sample at each corner of the domain. c) A combination of a) and b).}
\end{figure}\\ 
In this paper, the initial sample grid shown in figure \ref{fig:InitialSamples}c is employed, unless stated otherwise.\\

\noindent\textit{{Graph construction with neighbouring samples}}\\
\noindent The existing samples are connected based on Voronoi diagram construction \cite{barber_quickhull_1996}. A Voronoi diagram is a partitioning of a $d$-dimensional space into regions based on the distance to a specific set of samples \cite{goodman_handbook_1997}. Each Voronoi cell contains one sample, and the cell corresponds to all the points that are closer to this sample than to any other sample. The Voronoi diagram for the initial sample configuration is shown in figure \ref{fig:InitialSamplesVoronoi} (left).
\begin{figure}[!h]
\centering
\includegraphics[width = \textwidth]{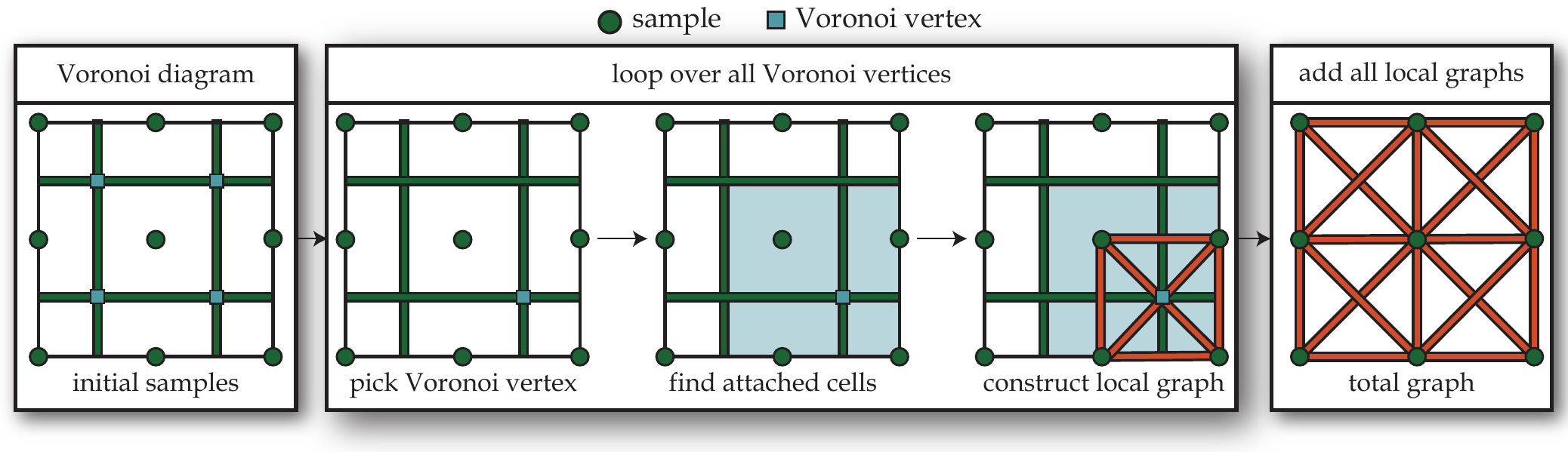}
\caption{\label{fig:InitialSamplesVoronoi} Algorithm for constructing the total graph of a set of sample points. (left) The Voronoi diagram with the corresponding Voronoi vertices. (middle) The local graphs connect the sample points of which the Voronoi cells have a common Voronoi vertex. (right) The total graph combines all the local graphs.}
\end{figure}\\
The boundary of each Voronoi cell contains several vertices. At each of these vertices the neighbouring Voronoi cells are determined. The midpoints of these neighbouring cells are then connected to obtain a local graph, see figure \ref{fig:InitialSamplesVoronoi} (middle). This local graph is similar to a Delaunay triangulation, but gives more isotropic behaviour on regular grids. The local graphs are subsequently connected together to obtain a global graph, see figure \ref{fig:InitialSamplesVoronoi} (right).\\

\noindent\textit{{Assignment of weights to edges}}\\
Not all the edges in the total graph are equally attractive for placing new samples. Edges that are long, and have large variation between QoI values and/or are in a region associated with a high PDF value, are good candidates for refinement. A weight function $w(\mbf{z}_1, \mbf{z}_2)$ assigns a weight to the edge between the two graph vertices $\mbf{z}_1$ and $\mbf{z}_2$. A low weight means that the edge is a good candidate for refinement and vice versa. Weighing based solely on either gradient \cite{jakeman_minimal_2013} or PDF \cite{witteveen_simplex_2010, witteveen_simplex_2012} is the most straightforward:
\begin{align}
w_{\text{PDF}}(\mbf{z}_1, \mbf{z}_2) &= \left(\rho\left(\frac{\mbf{z}_1 + \mbf{z}_2}{2}\right)||\mbf{z}_1 - \mbf{z}_2||_2\right)^{-1}\ , \label{eq:weight1}\\
w_{\text{grad}}(\mbf{z}_1, \mbf{z}_2) &= \left(\frac{|u(\mbf{z}_1) - u(\mbf{z}_2)|}{||\mbf{z}_1 - \mbf{z}_2||_2}||\mbf{z}_1 - \mbf{z}_2||_2\right)^{-1}\label{eq:weight2}\ ,
\end{align}
where $\rho$ is the PDF and $||\cdot||_2$ is the Euclidean distance. These weight functions are not always efficient, as they place most samples in either the smooth or the discontinuous regions. We therefore propose a new weight function that focuses on both the PDF and the gradient in the QoI:
\begin{equation}
w_{\text{PDF+grad}}(\mbf{z}_1, \mbf{z}_2) = \left(\left(\rho\left(\frac{\mbf{z}_1 + \mbf{z}_2}{2}\right) + \frac{|u(\mbf{z}_1) - u(\mbf{z}_2)|}{||\mbf{z}_1 - \mbf{z}_2||_2}\right)||\mbf{z}_1 - \mbf{z}_2||_2\right)^{-1}\ ,\label{eq:weight3}
\end{equation}
where both the PDF and gradient term contribute equally to the weight. We will compare this new weight function to the weight functions \eqref{eq:weight1} and \eqref{eq:weight2} in section \ref{sec:Results}. We multiply each weight function with the distance of the edge $||\mbf{z}_1 - \mbf{z}_2||_2$ to account for regions in the random space that have a low number of samples. After weighing all the edges, we normalise by dividing by the maximum weight. When using \eqref{eq:weight3}, the PDF and the gradient term are normalised separately, before adding them together and taking the reciprocal. The resulting weights are normalised again to ensure values in $[0, 1]$.\\

\noindent\textit{{Minimum spanning tree for refinement of sample grid}}\\
New samples are placed at the middle of the edges that have a sufficiently low weight. However, if all edges with sufficiently low weight are refined, undesirable clustering of samples may occur at early stages of the sampling procedure.

To prevent this, a Minimum Spanning Tree (MST) is used to obtain a subset of edges, such that this subset reaches all the samples, with a minimal total edge weight. The most important edges are contained in this MST, while still exploring a significant portion of the random space. The MST prevents the undesirable sample clustering at early stages. An edge in the MST is refined, if its weight is sufficiently low compared to the minimum weight among all edges, $w_{\text{min}}$:
\begin{equation}
w_i \leq c\ w_{\text{min}}\ ,
\label{eq:WeightCriterion}
\end{equation}
where $c>1$. The value of $c$ determines how many samples are added each iteration, i.e., low values of $c$ result in a low number of samples added and vice versa. Low values of $c$ produce the most accurate results, but many iterations are needed to reach a specified total number of samples. Iterations can become prohibitively expensive in high-dimensional random spaces. Therefore, in this paper we set $c=2$, which is a trade-off between the number of samples added and the number of iterations to be performed. The edge with minimum weight $w_{\text{min}}$ is not necessarily included in the MST, and will be added if it was not already included, to prevent the sampling algorithm from not adding any samples. \\ 

\noindent\textit{{Complete sampling algorithm}}\\
The complete sampling strategy (I) is an iterative procedure, which is illustrated in figure \ref{fig:SamplingStrategy}. The procedure starts with choosing the initial sample points. Next we loop over the three steps: \textit{graph construction}, \textit{edge weighing} and \textit{MST edge refinement}. The loop is terminated when the specified number of iterations $i_{\text{max}}$ has been performed or when the total number of sample points exceeds a specified threshold $N>N_{\text{max}}$.

The complete sampling strategy is illustrated in figure \ref{fig:SamplingStrategy}.
\begin{figure}[!h]
\centering
\includegraphics[width = \textwidth]{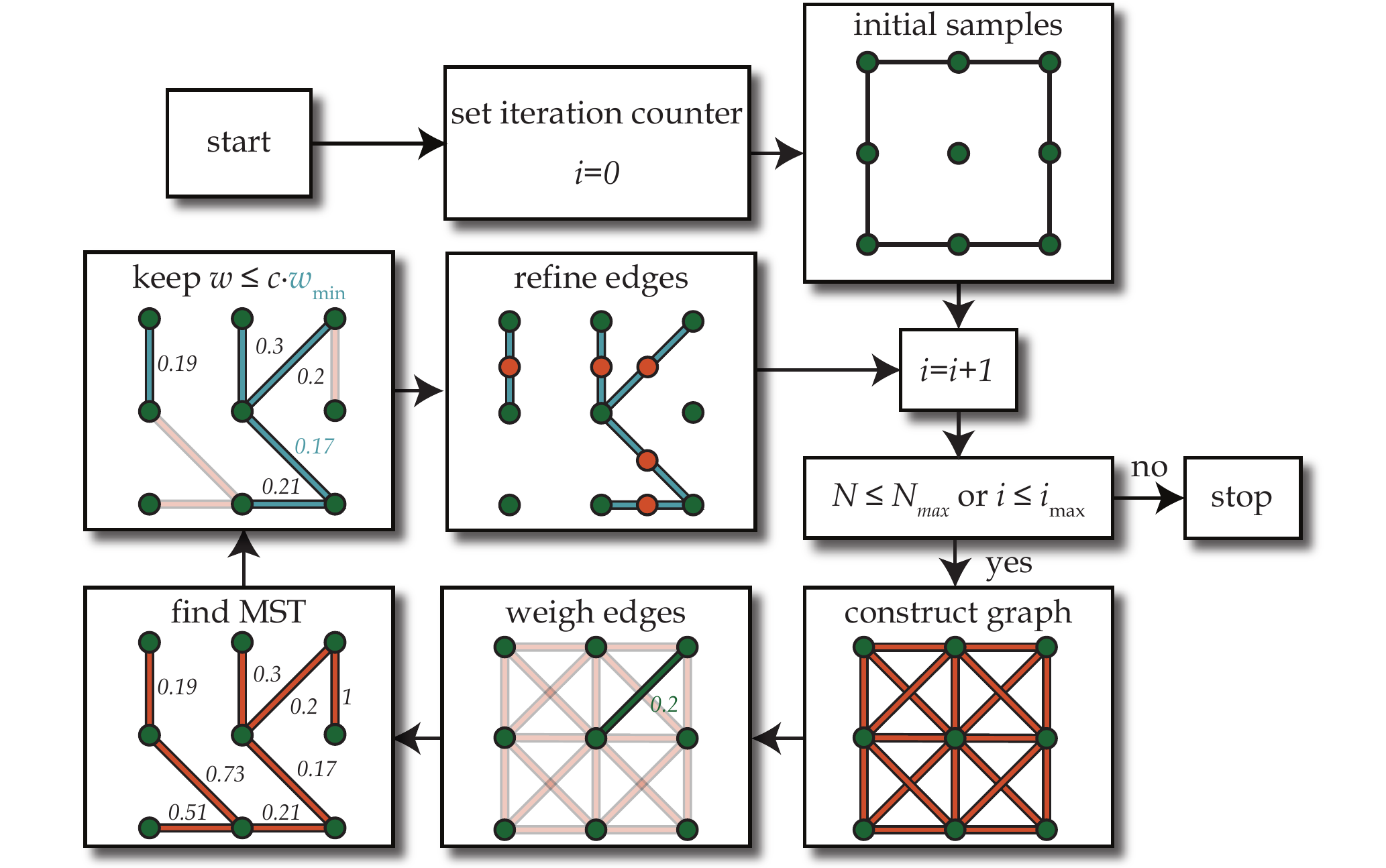}
\caption{\label{fig:SamplingStrategy} Schematic representation of the adaptive sampling strategy.}
\end{figure}
\subsection{Domain Decomposition}
\label{sec:DomainDecomposition}
\noindent The idea of the domain decomposition step in our method is to divide the random space into non-intersecting elements $E_i$, such that the sampled QoI values from section \ref{sec:SamplingProcedure} exhibit smooth behaviour locally in each element. The elements are constructed by first classifying the QoI values according to the QoI gradients. Second, the sample classes are separated by means of a classification boundary, based on SVMs. This classification boundary cuts the random space into several elements.\\

%

\noindent\textit{{Sample classification based on QoI gradients}}\\
\noindent An SVM determines a classification boundary based on a set of classified samples. Since a classification boundary is an approximation to a discontinuity in the QoI, the classification is based on the gradient between two neighbouring samples. A schematic representation of the classification procedure is shown in figure \ref{fig:ClassificationProcedure}. The choice for the threshold in the procedure is crucial. We have observed that the procedure from \cite{gorodetsky_efficient_2014} works very well and is therefore employed in this paper. It uses polynomial annihilation to estimate a jump value across an edge and labels two points the same class if the difference in function values is less than the jump value. An in-depth discussion of polynomial annihilation is discussed in \cite{jakeman_characterization_2011}.
\begin{figure}[!h]
\centering
\includegraphics[width = 0.9\textwidth]{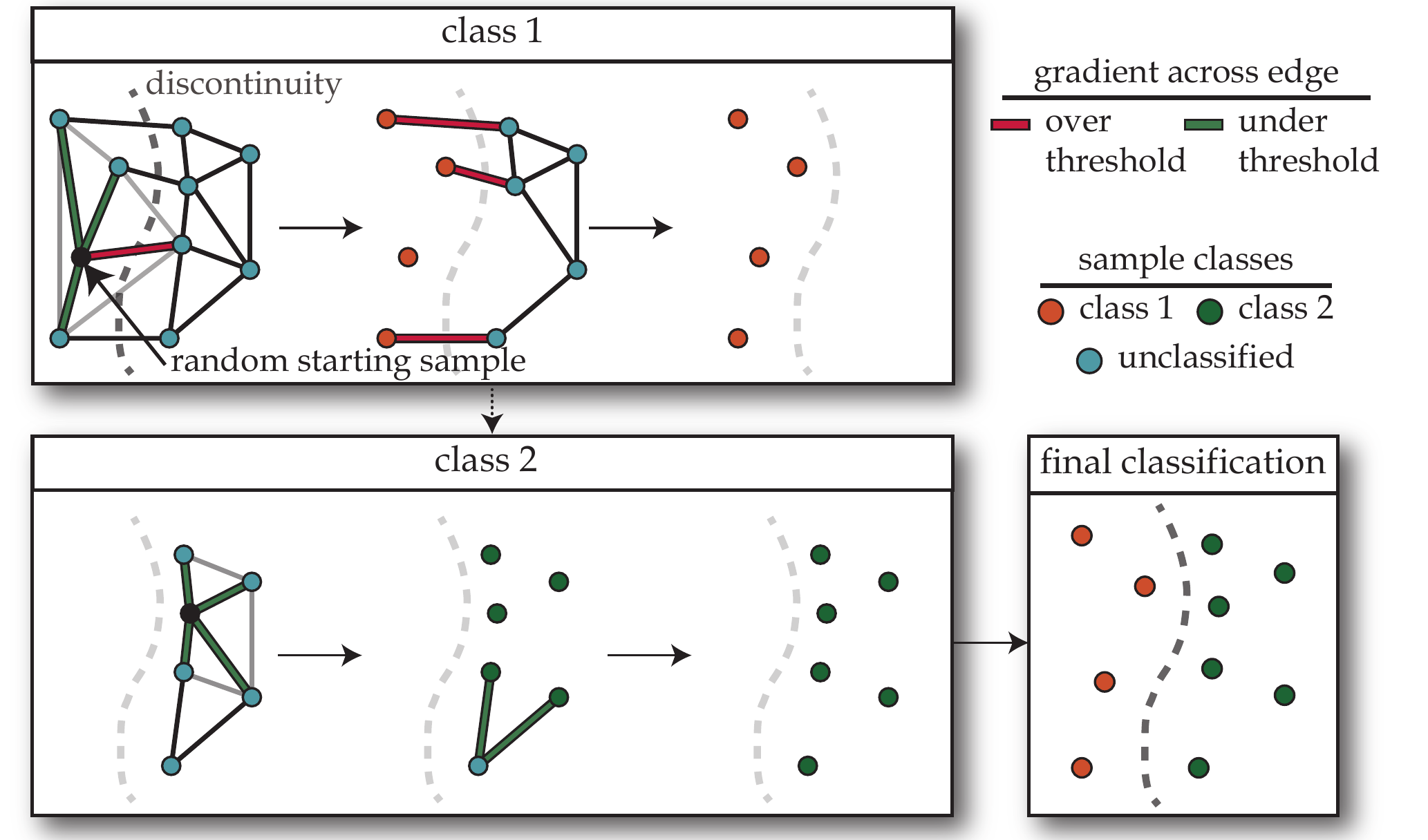}
\caption{\label{fig:ClassificationProcedure} Schematic representation of the classification procedure. The sample connectivity is given by the Voronoi construction performed on the final sampling grid from the adaptive sampling algorithm.}
\end{figure}\\

\noindent\textit{{Classification boundary from support vector machines}}\\
\noindent SVM is a supervised machine learning technique, used for building a classification boundary between samples that belong to different classes \cite{gorodetsky_efficient_2014, scholkopf_learning_2001}. SVM is used in this paper as it is defined by a convex optimisation problem, for which efficient methods are available \cite{chang_libsvm:_2011}, which makes it viable for high-dimensional problems.

Assume $N$ adaptive samples $\{\mbf{z}_i\}_{i=1}^n$ are classified into $N_c$ different classes $c_i$, where $c_i \in \{1,...N_c\}$ is the class belonging to sample point $\mbf{z}_i$. The idea behind an SVM is to construct a classifier $S_{\lambda}$ of the form:
\begin{equation}
S_{\lambda}(\mbf{z}) = \sum_{i=1}^{N_{sv}} \alpha_i K(\mbf{z}, \mbf{z}_i)\ ,
\label{eq:svmclassifier}
\end{equation}
where $\alpha_i$ is the coefficient associated with the sample point $\mbf{z}_i$, $\lambda$ a regularisation parameter, $N_{sv}$ the number of support vectors and $K$ a kernel \cite{burges_tutorial_1998}. If $\alpha_i > 0$, then $\mbf{z}_i$ is a support vector. Depending on the application, different kernels are available \cite{scholkopf_learning_2001}:
\begin{align}
K(\mbf{x}, \mbf{y}) &= <\mbf{x}, \mbf{y}>\ \ \text{(linear)}\label{eq:kernel1}\ ,\\
K(\mbf{x}, \mbf{y}) &= (\gamma<\mbf{x}, \mbf{y}> + c_t)^r\ \ \text{(polynomial)}\ ,\\
K(\mbf{x}, \mbf{y}) &= \exp(-\gamma ||\mbf{x} - \mbf{y}||_2^2)\ \ \text{(radial basis function)} \label{eq:rbfkernel}\ ,\\
K(\mbf{x}, \mbf{y}) &=\tanh(\gamma<\mbf{x}, \mbf{y}> + c_t)\ \ \text{(sigmoid)}\label{eq:kernel4}\ ,
\end{align}
where $<\cdot, \cdot>$ is the standard inner product in $\mathbb{R}^d$, $\gamma$ a regularisation constant, $c_t$ a translation constant and $r$ the polynomial degree. Radial basis function kernels \eqref{eq:rbfkernel} are commonly used and are also employed in this paper. The constant $\gamma$ is normally advised to be chosen as $1/{N_c}$ \cite{chang_libsvm:_2011}, but we will also investigate other choices in the result section of this paper. High values for $\gamma$ result in a classification boundary with a fine resolution, but this may result in overfitting. Low values for $\gamma$ result in a coarser estimation of the location of the classification boundary. The optimal value for $\gamma$ differs for different functions. In section \ref{sec:Results} an optimal value for $\gamma$ is found for a specific family of functions, which is important in our test-cases. Once $K$ has been chosen, the classifier $S_{\lambda}$ is obtained by solving the following least-squares problem:
\begin{equation}
S_{\lambda}(\mbf{z}) = \text{arg} \min_{S\in L_2(I_{\mbf{z}})} \left\lbrace \frac{1}{N}\sum_{i=1}^N \max(0, 1-c_i S(\mbf{z}_i)) + \lambda ||S||_{L^2(I_{\mbf{z}})}\right\rbrace\ .
\label{eq:svmproblem}
\end{equation}
The classification boundary is given by the $0$-contour of $S_{\lambda}$. It separates the different classes from each other with a hypersurface and is obtained with the LIBSVM library \cite{chang_libsvm:_2011}. The classification boundary decomposes the domain into several elements. Figure \ref{fig:SVMexample} shows an example of a classification boundary for two different classes. SVM can deal with multiple classes, and hence multiple discontinuities, quite easily, which makes it a suitable discontinuity finder for a wide range of QoIs. SVM in combination with the classification procedure is not yet able to detect discontinuities that run partially through the domain.

\begin{figure}[!h]
\centering
\includegraphics[width = \textwidth]{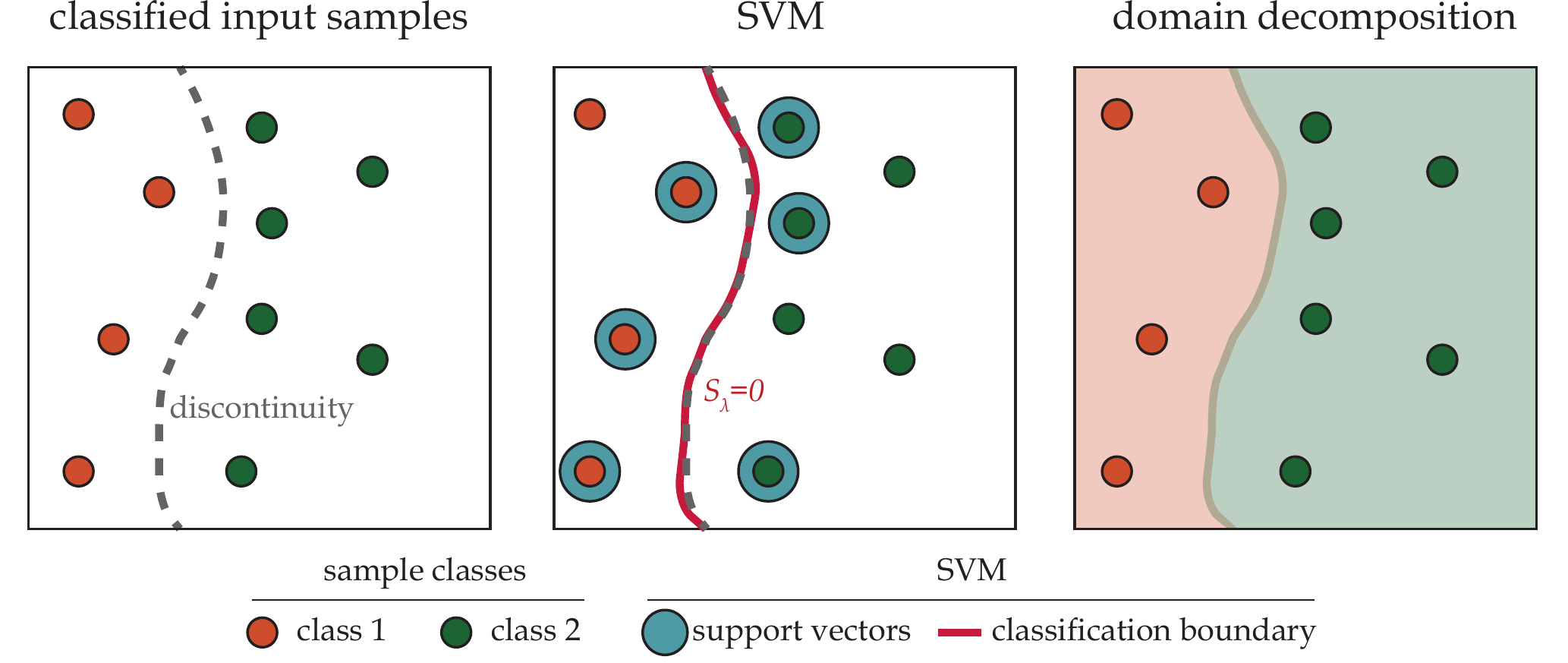}
\caption{\label{fig:SVMexample} An example of an SVM domain decomposition for two different classes.}
\end{figure}

\subsection{Local Approximations}
\label{sec:LocalApproximation}
\noindent The elements $E_i$ from the domain decomposition (\ref{sec:DomainDecomposition}) are arbitrarily shaped and the samples $\mbf{z}_i$ are distributed in such a way that interpolation is not a trivial task. Least orthogonal interpolation is able to perform interpolation on sample distributions on such arbitrarily shaped domains. For an in-depth discussion of least orthogonal interpolation, see \cite{narayan_stochastic_2012}. The sampling strategy (\ref{sec:SamplingProcedure}) does not necessarily choose points that are optimal for interpolation. Therefore, attempting to construct an interpolant on this set of interpolation nodes is not always a good idea and may produce unstable interpolants. We therefore use an extended version of the original least orthogonal interpolation, which selects a subset of samples that is better suited for stable interpolation \cite{jakeman_minimal_2013}. This enables the MST-ME method to place sample points in the random space, where we want to further resolve the QoI, without focusing on the stability of the interpolation.

We denote the least orthogonal interpolation operator by $I[\cdot]$, which operates on a subset of $(\mbf{z}_i, u(\mbf{z}_i))_{i=1}^n$, and we assume that the random space is decomposed into $N_c$ elements $E_i$. Each element $E_i$ comprises a single class $c_i$, which consists of the samples $(\mbf{Z}_i, \mbf{U}_i)$. The global approximation is given by, see equation \eqref{eq1:globalsolution}:
\begin{equation}
\tilde{u}(\mbf{z}) = \sum_{i=1}^{N_c} I[(\mbf{Z}_i, \mbf{U}_i)] \mathcal{I}_{E_i}(\mbf{z})\ ,
\label{eq1:globalsolution2}
\end{equation}
where $\mathcal{I}_{E_i}(\mbf{z})$ is the indicator function satisfying $\mathcal{I}_{E_i}(\mbf{z})=1$, if $\mbf{z}\in E_i$ and 0 otherwise.

\subsection{Complete Algorithm}
\noindent The complete MST-ME method, which comprises: sampling (I), domain decomposition (II) and local approximation construction (III), is shown in figure \ref{fig:CompleteAlgorithm}. Apart from the weight function and sampling threshold, no additional input from the user is required, which makes the method suitable for generic problems.
\begin{figure}[!h]
\centering
\includegraphics[width = \textwidth]{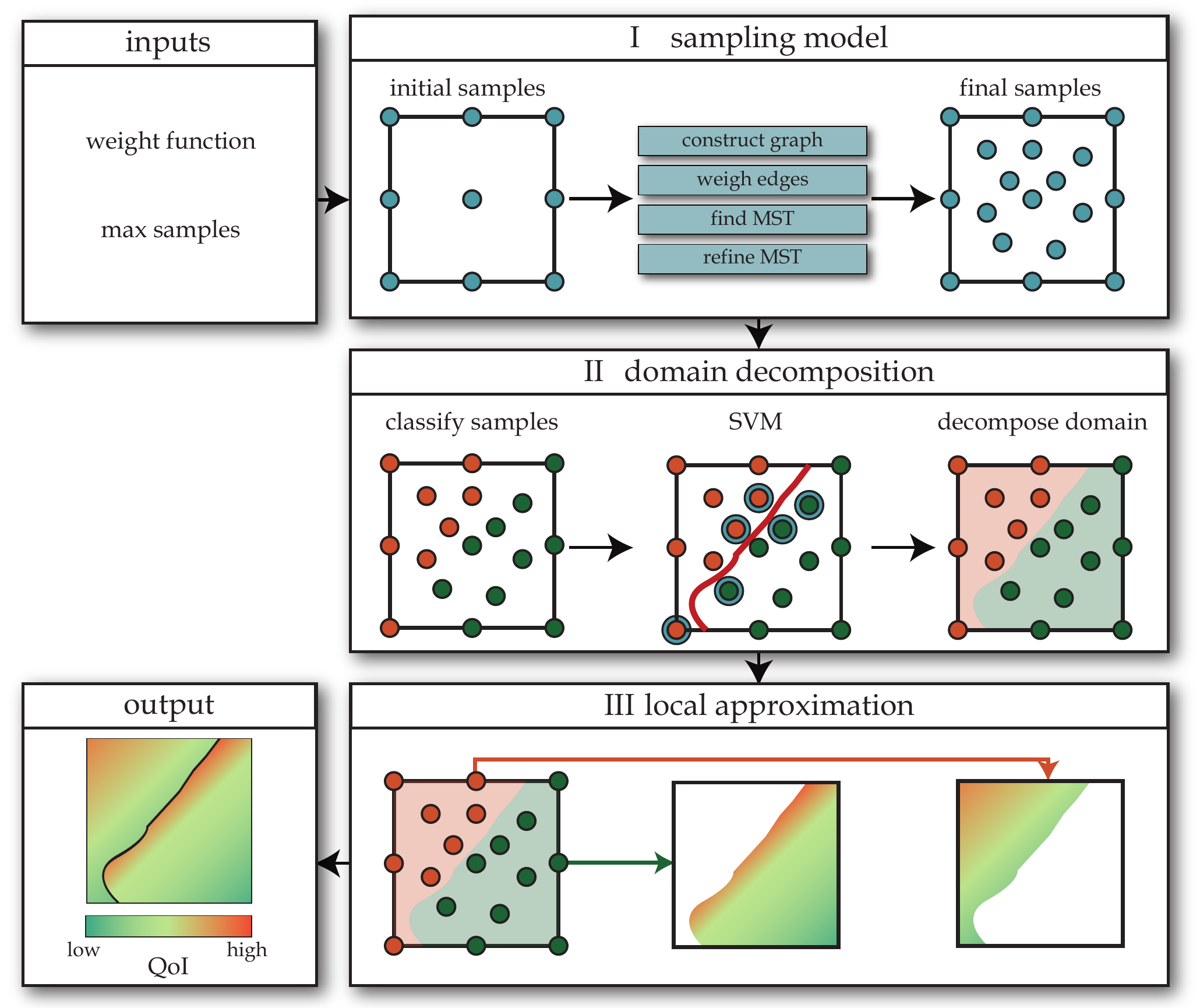}
\caption{\label{fig:CompleteAlgorithm} A schematic overview of the complete MST-ME method in a 2D random space.}
\end{figure}}\\

\noindent\textit{{Computational cost for $N$ samples in $d$-dimensional random space}}\\
\newcommand{\bigO}{\ensuremath{\mathcal{O}}}
\noindent The computational cost of the sampling strategy (I) is determined by the cost of computing the Voronoi diagrams and finding the MST. Computing a Voronoi diagram on $N$ samples in $\mathbb{R}^d$ can be done in $\bigO(N\log(N) + N^{\lceil d/2\rceil})$ time \cite{fortune_sweepline_1986}. In the worst-case scenario, only a single sample is added in each iteration, and hence we have to compute $N$ Voronoi diagrams on $N$ samples. The maximum computational cost is therefore $\bigO(N^2\log(N) + N^{\lceil d/2\rceil + 1})$. Computing the MST can be done with Prim's algorithm \cite{prim_shortest_1957}, which has an algorithmic complexity of $\bigO(|E|\log(|V|))$, where $|E|$ is the total number of edges and $|V|$ the total number of samples. An upper bound for $|E|$ and $|V|$ is given by $N(N+1)/2$ and $N$, respectively. Again, in the worst case scenario we have to perform $N$ iterations, which can be done in $\bigO(N^3\log(N))$ time.

The computational cost of the domain decomposition (II) is based on the complexity of the classification procedure and the SVM. Classification of $N$ samples can be performed in $\bigO(N^2)$ time. The SVM has a complexity ranging between $\bigO(N^2)$ and $\bigO(N^3)$, depending on the number of classes and the kernel \cite{chang_libsvm:_2011}. Hence, domain decomposition can be performed in $\bigO(N^3)$ time and is independent of $d$.

The local approximation (III) uses LU and QR-decomposition for determining the interpolant. Computing the QR-decomposition is the dominant factor in (III), and it can be done in $\bigO(N^3)$ time, using the standard implementation in Matlab. If least orthogonal interpolation is employed in an adaptive fashion, the complexity increases to $\bigO(N^4)$, as $N$ QR-decompositions are computed in the worst case.

The complexity of the complete algorithm is determined by the complexities of (I), (II) and (III). This results in an overall complexity of $\bigO(N^4)$ if $d\leq 5$ and $\bigO(N^3\log(N) + N^{\lceil d/2\rceil + 1})$ otherwise.\\

%
\section{Results}
\label{sec:Results}
\noindent In this section we present multiple examples that illustrate the robustness and flexibility of the MST-ME method. For computing the error between the exact surrogate and approximation we use the following weighted $L_{2,\rho}$-norm:
\begin{align}
\norm{\tilde{u} - u}_{2,\rho}^2 &= \frac{1}{N_{MC}}\sum_{i=1}^{N_{MC}} \rho(\mbf{z}^{MC}_i)\cdot|\tilde{u}(\mbf{z}^{MC}_i) - u(\mbf{z}^{MC}_i)|^2\ ,
\label{eq:SWEL2}
\end{align}
where the surrogate model $\tilde{u}$ is constructed using MST-ME and evaluated at Monte Carlo sample locations $\mbf{z}_i^{MC}$ in the random space. The exact solution $u$ is the evaluation of the model sampled at the same Monte Carlo samples. We multiply the difference between the surrogate and exact solution with the PDF $\rho$ to emphasise on regions in the random space that are likely to occur. Samples that are within a distance from the discontinuity, which is lower than the minimum distance of the adaptive MST-ME samples, are discarded. Discarding samples is motivated by the fact that these points lie below the resolution of the SVM discontinuity detection, where the fidelity of the classification is questionable \cite{jakeman_minimal_2013}.

The first example shows the approximation of three 2D, piecewise constant functions. This example focuses solely on the domain classification and shows convergence of the SVM domain decomposition (II). A second example shows the approximation of 1D and multidimensional Genz functions \cite{genz_testing_1984}. The third test-case shows the MST-ME method applied to a more complicated model, which is defined by an underlying PDE, namely, the shallow water equations. Lastly, we apply the MST-ME method to a 3D dam break problem governed by incompressible fluid-flow equations, to show that our methodology can be used for a complex engineering problem.
\subsection{Domain Classification}
\noindent The accuracy of the SVM domain decomposition (II) is investigated as a function of the parameter $\gamma$. \\

\noindent\textit{{SVM domain decomposition works best for discontinuities without corners}}\\
The SVM domain decomposition (II) is tested by approximating three piecewise constant 2D functions. The functions have a discontinuity in the shape of a circle, a rectangle and a triangle, respectively. The sampling algorithm (I), with a weight function that focuses solely on the gradient \eqref{eq:weight2}, is used to determine the sample locations. The SVM domain decomposition uses the radial basis function kernel \eqref{eq:rbfkernel} in which $\gamma$ is set to the advised value $1/N_c$ \cite{chang_libsvm:_2011}, where $N_c$ is the number of classes. The sample locations for the circle test-case are shown in figure \ref{fig:ResultsDC1} (left), along with the definition of the error measure (the misclassified portion). Clustering of samples appears around the circular-shaped discontinuity location, because the gradient based weight function \eqref{eq:weight2} does not lead to refinement if there is no intersection with the discontinuity. The sample locations show symmetry, but after 5 iterations, the symmetry is lost slightly, although this is not noticeable in figure \ref{fig:ResultsDC1} (left).
\begin{figure}[!h]
\centering
\includegraphics[width = \textwidth]{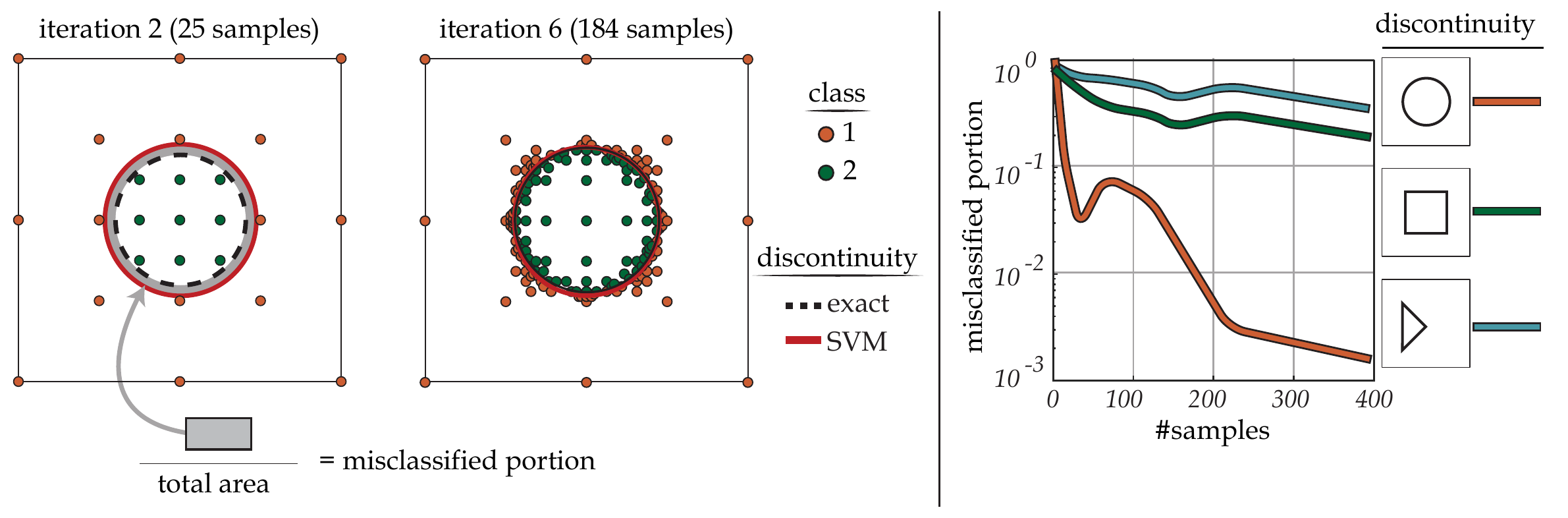}
\caption{\label{fig:ResultsDC1} (left) Sample locations at different iterations for the circular test-case, with $\gamma = \frac{1}{2}$. (right) Convergence of the misclassified portion with an increasing number of samples.}
\end{figure}\\
The domain decomposition error as a function of the number of samples is shown in figure \ref{fig:ResultsDC1} (right).

The misclassified portion of the domain decreases rapidly with the increasing number of samples, but is in general not monotonically decreasing. Discontinuity lines with sharp corners, in particular the square and triangle, are in general harder to approximate for the SVM discontinuity finding, than smooth discontinuity lines. The value for $\gamma$ may significantly influence the accuracy of the discontinuity approximation, and the advised value $1/N_c$ is in general not optimal. Therefore we search for a value of $\gamma$ that is optimal for the remaining test-cases.\\

\noindent\textit{{The optimal value for $\gamma$ for our test-cases is $3/N_c$}}\\
We now investigate the optimal value for $\gamma$ by testing several candidate values for a large set of piecewise constant functions. The parameter $\gamma$ may influence the accuracy of the SVM domain classification. Many discontinuous QoIs in engineering applications possess a discontinuity without corners, i.e., the exact discontinuity surface is a smooth hypersurface. For such discontinuities, the value of $\frac{1}{N_c}$ might not be optimal and therefore we search for an optimal value among a set of candidates that are perturbations on the advised value: $\{ \frac{1}{N_c+4}, \frac{1}{N_c+3}, \frac{1}{N_c+2}, \frac{1}{N_c+1}, \frac{1}{N_c}, \frac{2}{N_c}, \frac{3}{N_c}, \frac{4}{N_c}, \frac{5}{N_c}\}$. A set of $10^6$ piecewise constant functions, possessing up to 3 discontinuities, on the domain $[-1,  1]^2$, is randomly generated. The discontinuity location is given by up to 3 non-intersecting polynomial parametric curves up to degree 5, which have random coefficients. For each of these functions, an SVM domain decomposition is performed with each of the possible candidates for $\gamma$. This domain decomposition is performed on 50 adaptively sampled points from the MST-ME method, based on the weight function \eqref{eq:weight2} and the initial configuration shown in figure \ref{fig:InitialSamples}c. We add all correctly classified portions for each $\gamma$-candidate, and divide by $10^6$ to obtain an average correctly classified portion for all the randomly generated functions. The number of adaptive samples influences the average correctly classified portion, but a similar trend between the $\gamma$-candidates is obtained. The results are shown in figure \ref{fig:ResultsDC3}.
\begin{figure}[!h]
\centering
\includegraphics[width = \textwidth]{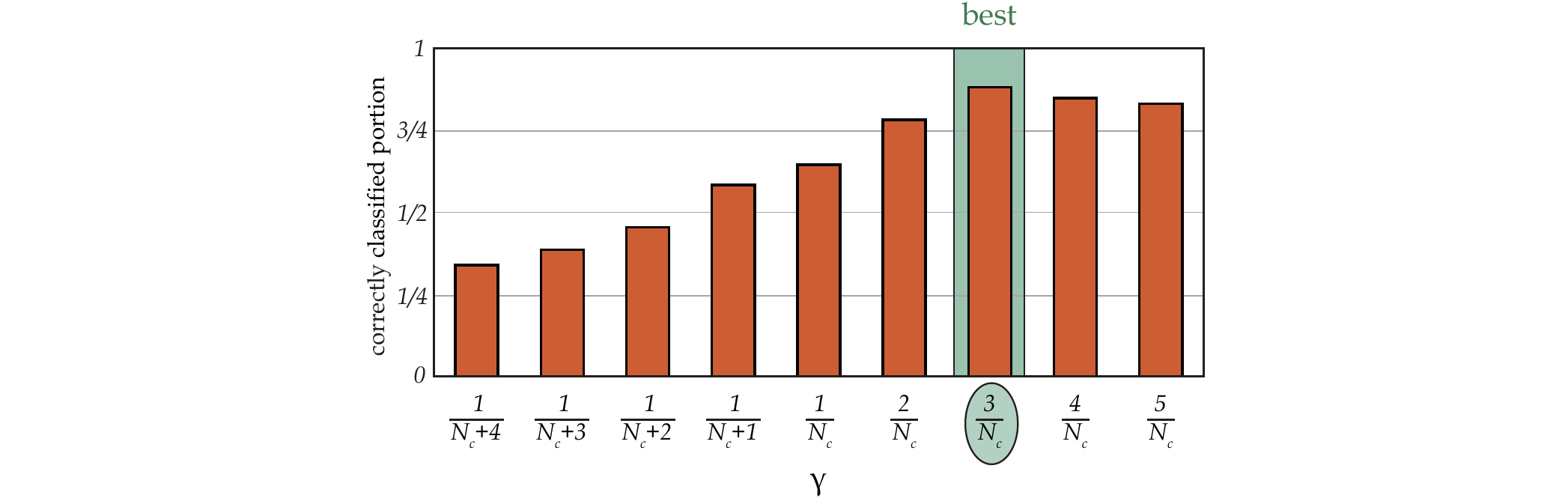}
\caption{\label{fig:ResultsDC3} Average correctly classified portion for different $\gamma$-values.}
\end{figure}\\
Figure \ref{fig:ResultsDC3} shows that $\gamma=\frac{3}{N_c}$ is the most accurate choice for the generated family of piecewise constant functions. The inconsistency with the value for $\gamma$ suggested in literature $(\frac{1}{N_c})$ and our value is possibly due to the fact that we limit ourselves to a family of discontinuous functions, which possess no sharp corners in the discontinuity surface. Hence, the value $\frac{3}{N_c}$ might not be the optimal value for other families of discontinuous functions. The discontinuities considered in the remainder of this paper fall within this category and $\gamma=3/N_c$ is therefore used in the remainder of this paper.
\subsection{Genz Functions Approximation}
\noindent To illustrate the accuracy/efficiency, the proposed MST-ME method is applied to a standard benchmark problem, namely, approximation of 1D Genz functions \cite{genz_testing_1984}.\\

\noindent\textit{{Edge weighing based on PDF and gradient is most robust}}\\
We consider the following Genz functions:
\begin{align}
g_1(x, \alpha) &= \cos(\alpha x)\ ,\\
g_3(x, \alpha, \beta) &= \left(\frac{1}{1+\alpha x}\right)^{1+\beta}\ ,\\
g_5(x, \alpha, \beta) &= \exp(-(\alpha|x|-\beta))\ ,\\
g_6(x, \alpha, \beta) &= \left\{\begin{array}{l l l}
0, &x>\beta\ ,\\
\exp(\alpha x), & \text{otherwise}\ .
\end{array}\right.
\end{align}
A uniform PDF is assumed on the interval $[-1, 1]$ and the initial grid consists of the two end points plus the middle point of this interval. As a reference, the MST-ME solution is compared with the stochastic collocation (SC) solution on a Gauss-Legendre grid.
\begin{figure}[!h]
\centering
\includegraphics[width = \textwidth]{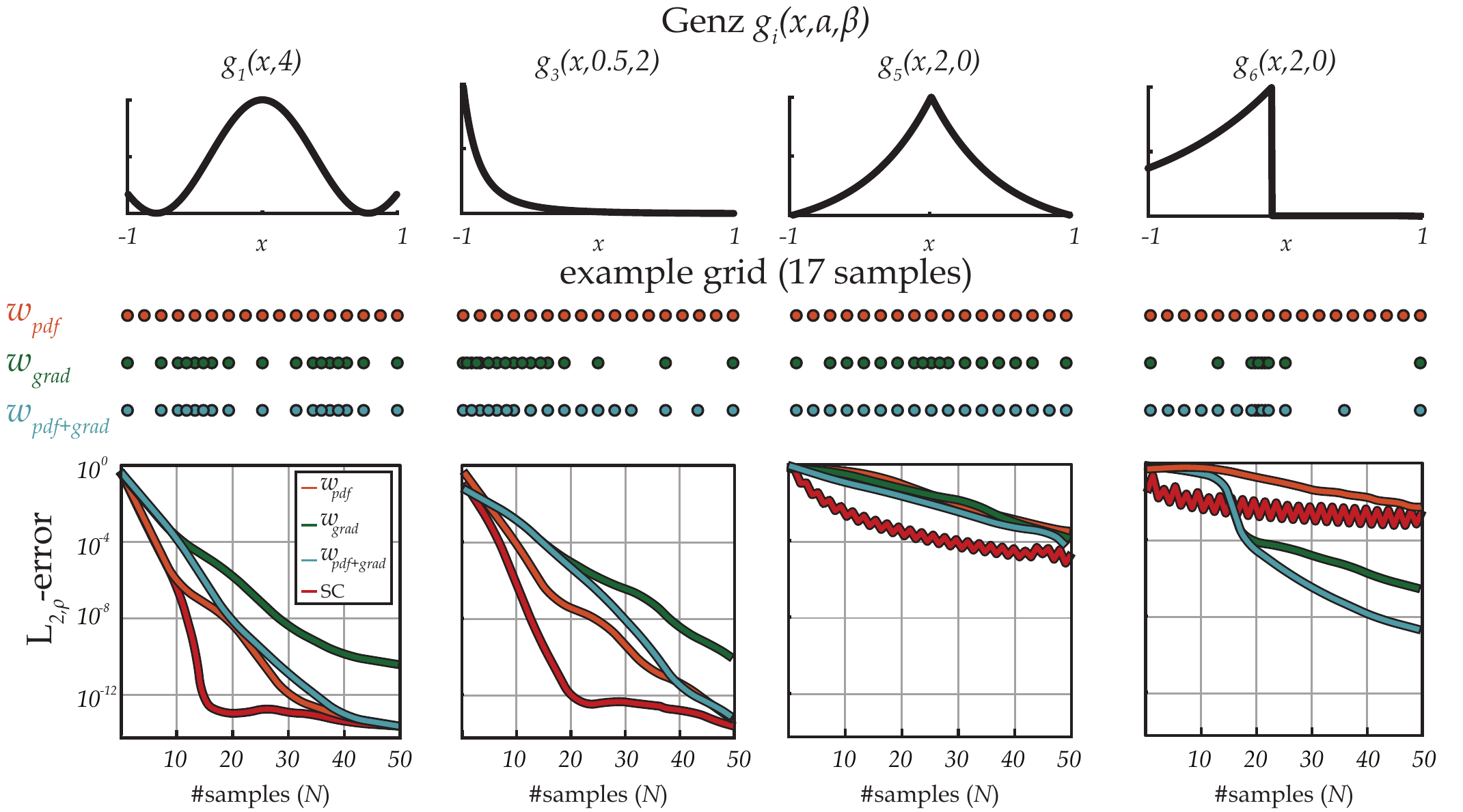}
\caption{\label{fig:ResultsGENZ1} Error of the approximation with respect to the exact function.}
\end{figure}\\

Figure \ref{fig:ResultsGENZ1} shows that weighing based on the PDF only results in a uniformly spaced sample grid. Interpolation on such a grid is in general not a good idea, as it may produce unstable interpolants. The least orthogonal interpolation method partially circumvents this issue by choosing a subset of these samples in constructing an interpolant. Consequently, the smooth Genz functions $g_1$ and $g_3$ are well approximated, but $g_5$ and $g_6$ are not. Weighing based on the gradient alone leads to improved results for the discontinuous function $g_6$, but leads to less accurate results for the smooth functions $g_1$ and $g_3$. The standard stochastic collocation method performs well in smooth cases, but converges slowly and also shows an oscillating error in some cases, due to the Gauss-Legendre grid that includes the middle point of the domain only for an odd number of samples. In contrast, the weight function based on both the PDF and the gradient performs the best overall, by keeping track of the discontinuity location, while still maintaining a sample distribution which resolves parts of the random space away from the discontinuity. In the remainder of this paper we therefore use weighing based on both PDF and gradient, equation \eqref{eq:weight3}.

Notice that all weight functions show slow convergence for approximating $g_5$. This is due to the absence of the second derivative in the weight function \eqref{eq:weight3}. By basing the classification on the second derivative in the QoI, we can circumvent the slow error convergence in presence of a discontinuity in the first derivative. However, we will not have any discontinuities in the first derivatives for the remaining test-cases in this paper. We therefore use classification based on the first derivative only in the remainder of this paper.\\

\noindent\textit{{The underlying PDF changes the sample grid}}\\
The effect of the underlying PDF is now investigated. The two PDFs that we consider are a symmetric and an asymmetric $\beta$-distribution, with parameters $(10, 10)$ and $(2, 7)$, respectively. The support of both PDFs is scaled to $[-1, 1]$ and we use the uniform distribution as a reference. The error convergence for $g_1$ and $g_6$ is shown in figure \ref{fig:ResultsGENZpdf}.
\begin{figure}[!h]
\centering
\includegraphics[width = \textwidth]{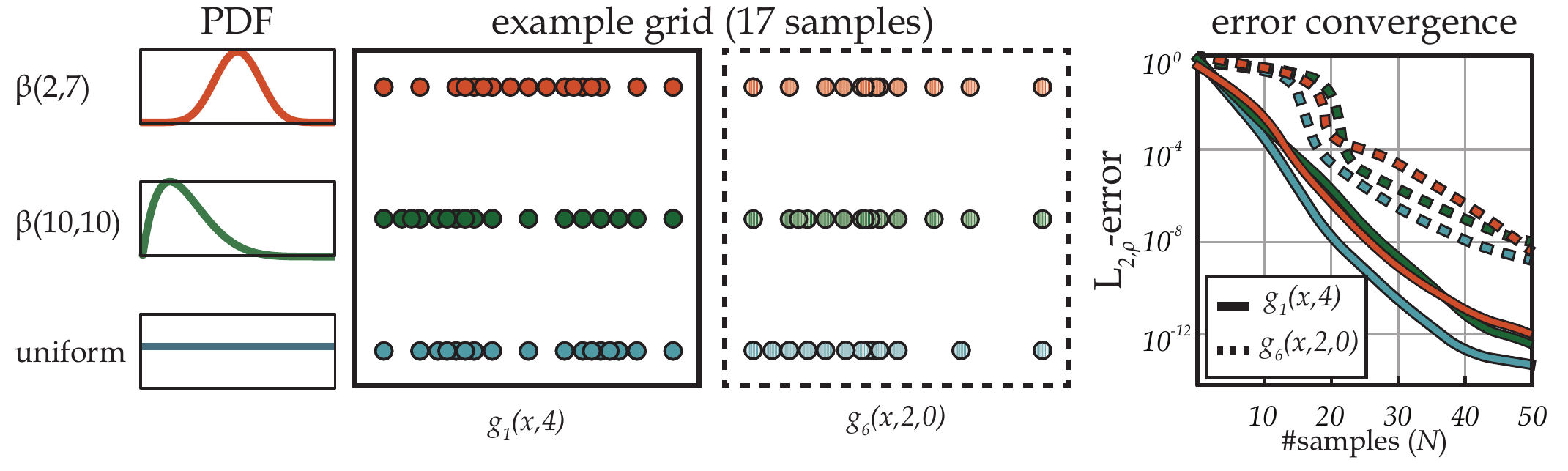}
\caption{\label{fig:ResultsGENZpdf} Error convergence for different choices of PDFs.}
\end{figure}\\
The error convergence is similar to the error convergence for the uniform distribution. Again the sample grid is not ideal for interpolation, but the adaptive least orthogonal interpolation circumvents this by using a subset of nodes.\\

\noindent\textit{{Samples per dimension is constant}}\\
To investigate how the MST-ME method performs in multiple dimensions, the error is plotted for a smooth Genz function $g_1$ and a discontinuous Genz function $g_6$, with increasing dimension $d$, in figure \ref{fig:ResultsGENZ2} (left). These multidimensional Genz functions are tensor products of the 1D Genz functions. The error for each dimension is based on 1000 adaptively sampled points with an initial sampling configuration equal to the one shown in figure \ref{fig:InitialSamples}c. The number of required samples needed to attain a specific accuracy is also plotted as a function of the dimension of the random space (right).
\begin{figure}[!h]
\centering
\includegraphics[width = \textwidth]{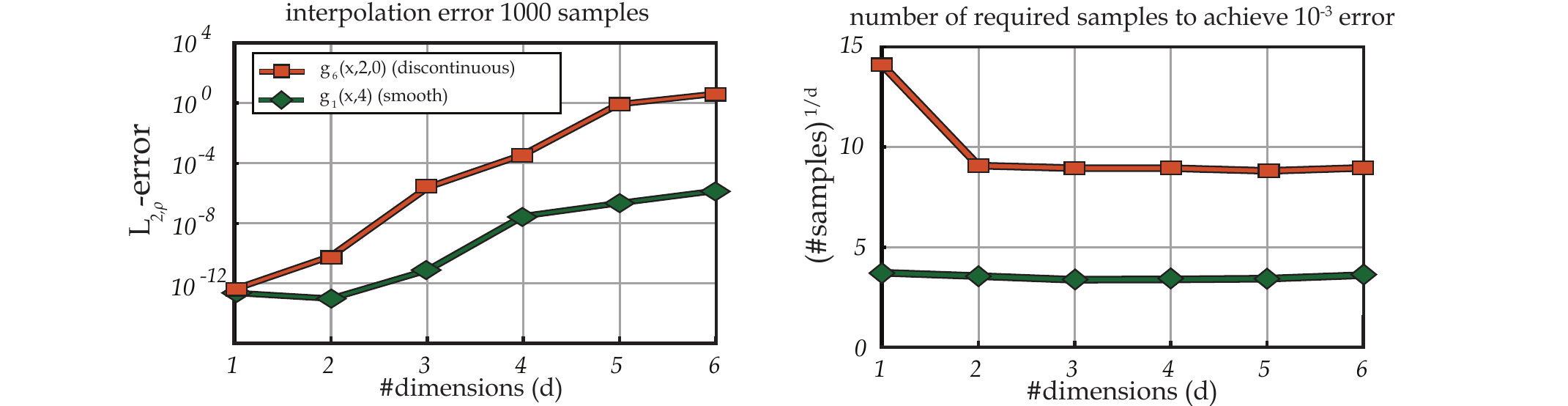}
\caption{\label{fig:ResultsGENZ2} (left) Error of the approximation with respect to the exact function as a function of the dimension. (right) The number of samples needed to attain a certain accuracy as a function of the dimension.}
\end{figure}\\
Apart from the first point in the discontinuous case, the average number of samples per dimension is approximately constant, which is the well-known curse of dimensionality. Since there is no strong change in the number of samples per dimension, the MST-ME method is applicable to problems with a low to medium number of uncertain inputs.
\subsection{Shallow Water dam break}
\noindent We study the performance of the MST-ME method applied to a system of 1D conservation laws. This system consists of the 1D shallow water equations (SWEs), which describe the inviscid flow of a layer of fluid with free surface, under the action of gravity, with the thickness of the fluid layer small compared to the other length scales \cite{vreugdenhil_numerical_2013}:
\begin{equation}
\der{}{t}\vecbrace{h \\hv} + \der{}{x}\vecbrace{hv \\hv^2 + gh^2/2} = \mbf{0}\ ,
\label{eq:SWE1D}
\end{equation}
where $h$ is the free surface height (thickness of the fluid layer), $v$ the velocity, and $g$ the acceleration of gravity. The initial condition for the system of PDEs is given by:
\begin{equation}
\vecbrace{h \\ v}(x, t=0) = \left\lbrace \begin{array}{ll}
\vecbrace{h_l \\ v_l},& x\leq 0\ ,\\
\vecbrace{1 \\ 0},& x> 0\ ,\\
\end{array} \right. 
\end{equation}
leading to a Riemann problem shown in figure \ref{fig:ResultsSWESchematic}. The solution of the Riemann problem for these initial conditions can be computed exactly when working on an infinite spatial domain \cite{leveque_finite_2002}. The solution consists of two characteristic waves travelling through the spatial domain, see figure \ref{fig:ResultsSWESchematic}. Each wave is a shock or rarefaction wave. When solid boundary conditions are imposed at $x=\pm 1$, an exact solution cannot be obtained for all initial solutions. We therefore employ a finite volume method with an exact Riemann solver \cite{toro_riemann_2013} to compute the cell face fluxes, and solve the SWEs using 256 cells. A Crank-Nicolson scheme is used to integrate the SWEs in time and a ghost-cell method with reflective properties is used for the boundaries. The initial left state $(h_l, v_l)$ is assumed to be uncertain and uniformly distributed between $[1.5, 2.5]$ and $[-0.5, 0.5]$, respectively, i.e.:
\begin{equation}
\mbf{z} = \vecbrace{h_l \\ v_l} \sim \vecbrace{\mathcal{U}( 1.5, 2.5 ) \\ \mathcal{U}( -0.5, 0.5 )} \ .
\end{equation}
The uncertainty in the initial conditions is large in order to ensure that we get different characteristic behaviour of the QoI. The average thickness of the fluid layer is not shallow compared to the domain length, but this is not important for showing the performance of the MST-ME method. The QoI $u$ is defined as the fluid height at $x=-1$ at a certain time $t^*$. A schematic representation of this set-up is shown in figure \ref{fig:ResultsSWESchematic}.
\begin{figure}[!h]
\centering
\includegraphics[width = \textwidth]{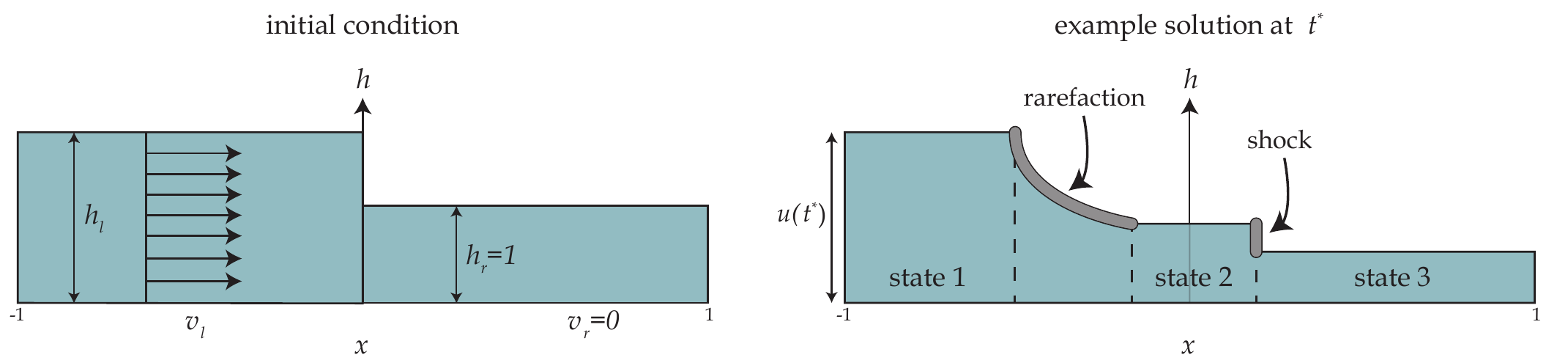}
\caption{\label{fig:ResultsSWESchematic} Schematic representation of the shallow water test-case.}
\end{figure}\\
Notice that the QoI is time dependent and the characteristics of this QoI will change significantly as time progresses. Either a transition between a shock and rarefaction wave, or a difference in wave speeds can result in a discontinuity in $u$. This allows us to study the robustness of the MST-ME method, as this test-case comprises both smooth and discontinuous QoI responses. We emphasise that the constructed surrogate for the QoI at a $t^*$ cannot be reused for other time instances, because the MST-ME method uses the QoI at the current time $t^*$ as a measure to place new samples.\\

\noindent\textit{{MST-ME detects if a function is smooth or discontinuous automatically}}\\
The MST-ME method is used for three different QoIs, $u=h(x=-1, t^*\in\{1.67, 4.16, 2.21\})$, which correspond to a mildly non-linear, highly non-linear and close to discontinuous QoI, respectively. The surrogate model and sample grids after 10 iterations are shown in figure \ref{fig:ResultsSWE1}.
\begin{figure}[!h]
\centering
\includegraphics[width = \textwidth]{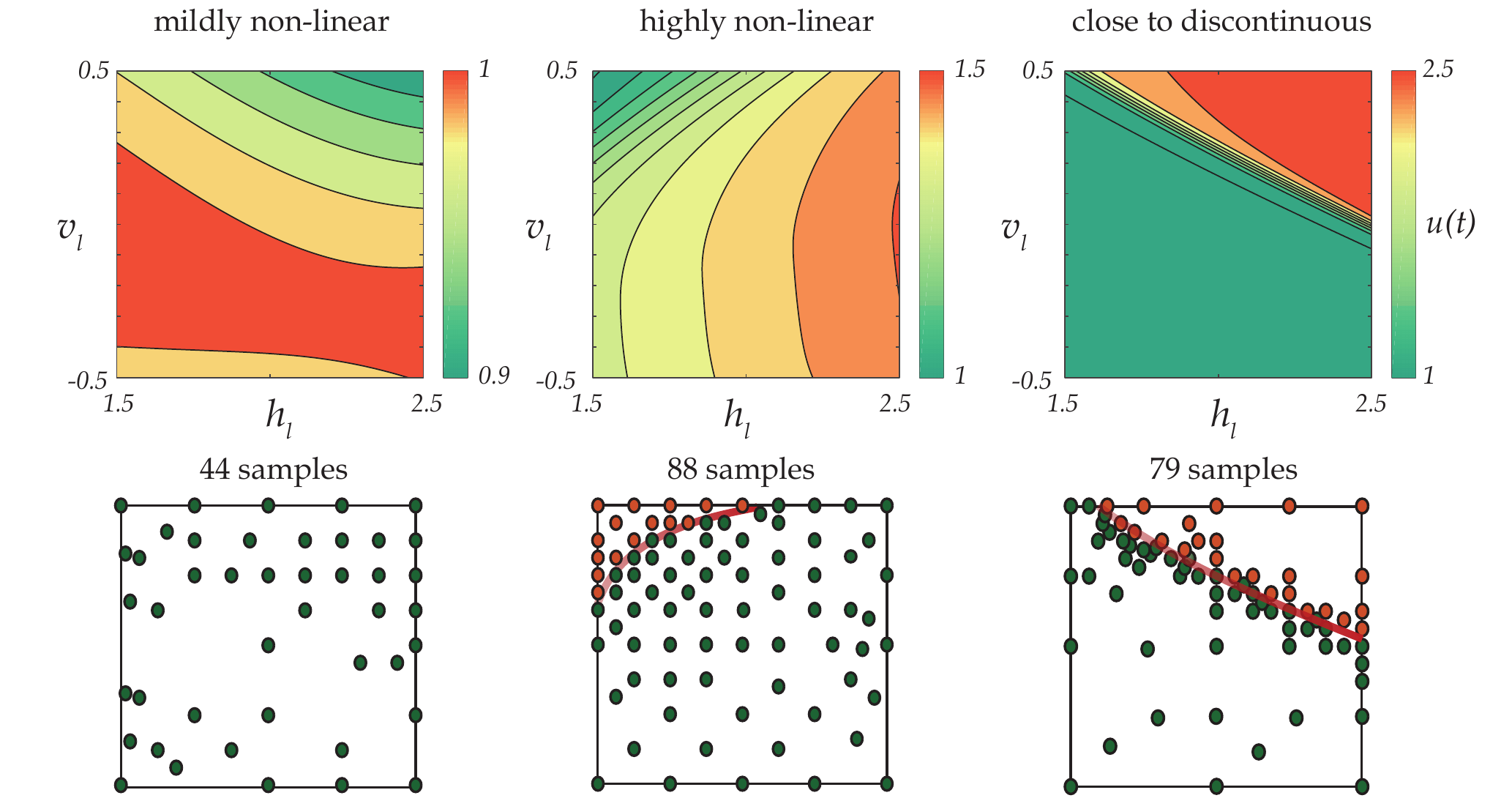}
\caption{\label{fig:ResultsSWE1} The three QoI surrogate models and corresponding sample grids after 10 iterations.}
\end{figure}\\
The discontinuity in the QoI at time $t^*=2.21$ is caused by a shock wave, which hits the left boundary for certain values in the random space, but does not yet hit the left boundary for other values in the random space. 

To investigate the accuracy of the MST-ME method, we determine the convergence. The error is based on $10^6$ Monte Carlo samples. The convergence of the $L_{2,\rho}$-error \eqref{eq:SWEL2} is shown in figure \ref{fig:ResultsSWE2}. 
\begin{figure}[!h]
\centering
\includegraphics[width = \textwidth]{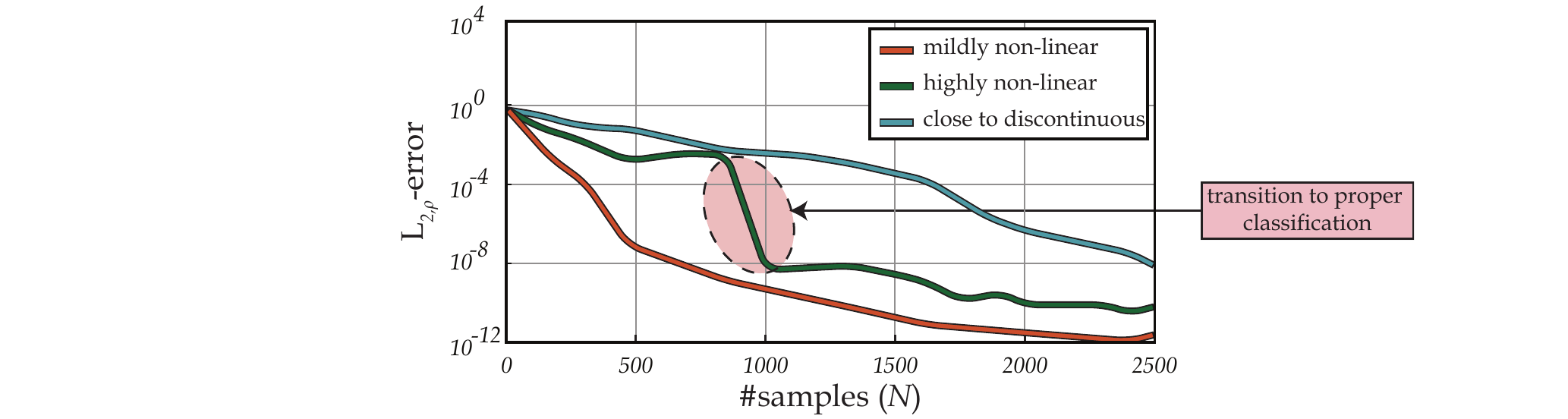}
\caption{\label{fig:ResultsSWE2} Convergence of the MST-ME solution in the $L_{2,\rho}$-error \eqref{eq:SWEL2}.}
\end{figure}\\
The results show that the error as a function of the samples decays fast for the mildly and highly non-linear case, as expected. The highly non-linear case shows a sudden drop in the error, which is caused by transition in the domain decomposition. First the classification procedure detects a large enough jump in the sampled QoI to conclude that there is a discontinuity present in the QoI. As the MST-ME progresses, samples are added in the area of the possible discontinuity, until the jump in the QoI values becomes small enough to classify the samples properly. This transition to a correct classification explains the sudden drop in the error for the highly non-linear case. MST-ME automatically detects the smoothness of the QoI, as the number of samples increases.
\subsection{3D dam break}
\noindent As a last test-case, the MST-ME method is applied to a complex engineering problem, namely a 3D fluid dam break problem with parametric initial conditions. This test-case is similar to the shallow water dam break, with the difference that it describes fluid motion in 3D and contains more physics, i.e., viscous effects and no shallow-water assumption. Dam break problems are commonly used as benchmark in for example sloshing applications \cite{hopfinger_liquid_nodate} and have been extensively studied \cite{hopfinger_liquid_nodate,marsooli_3-d_2014,larocque_3-d_2013}. MST-ME can be used to gain physical insight for this parametric problem, by constructing a surrogate model in the full parameter space, and this is our goal in this test-case. We do not focus on convergence, as it is computationally infeasible to construct a reliable reference solution. A schematic representation of the test-case is shown in figure \ref{fig:Results3DSchematic}.
\begin{figure}[!h]
\centering
\includegraphics[width = 0.9\textwidth]{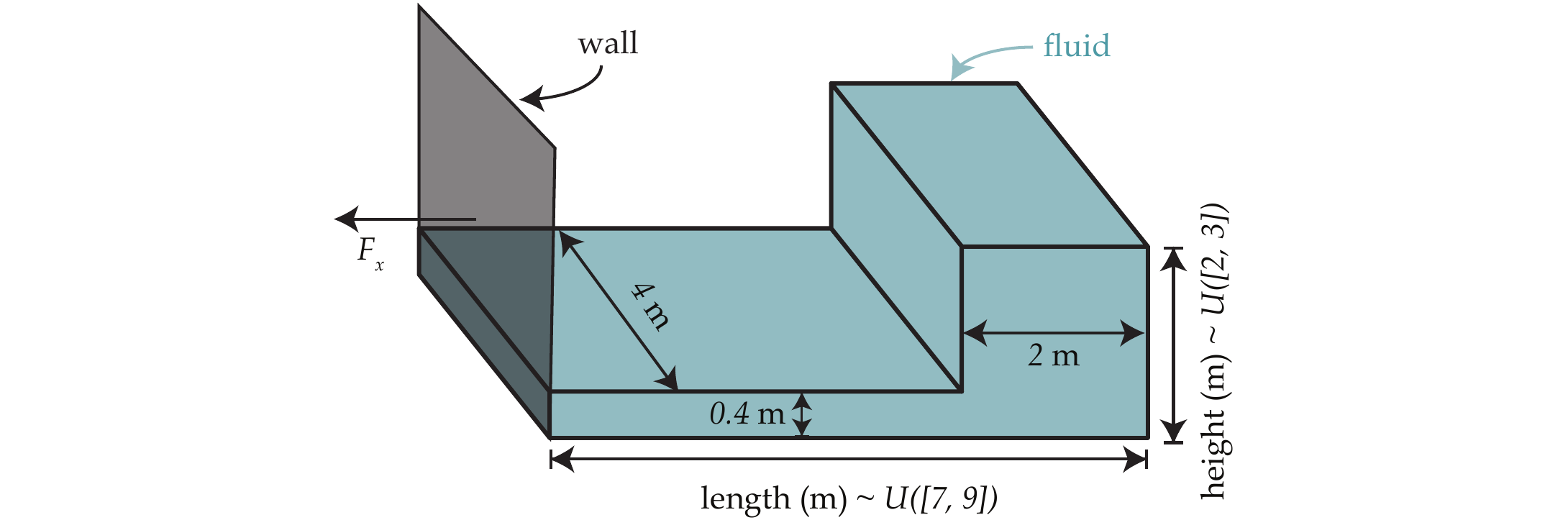}
\caption{\label{fig:Results3DSchematic} Schematic of 3D dam break.}
\end{figure}\\
We consider two uncertain parameters, namely the length of the tank and the height of the right column of liquid. When the fluid is released, it starts flowing to the left side of the domain and impacts the wall, shown in figure \ref{fig:Results3DSchematic}. The QoI is the maximum perpendicular wall force component $F_x$ during first impact. Smoothed particle hydrodynamics (SPH) is used to simulate the free surface flow, induced by the initial condition in figure \ref{fig:Results3DSchematic}. We use an open-source SPH solver, DualSPHysics \cite{crespo_dualsphysics:_2015}. The acceleration of gravity is set to $9.81\text{m}/\text{s}^2$ and the kinematic viscosity $\nu$ is set to $10^{-6}\text{m}^2/\text{s}$, which is the value for water at room temperature. Surface tension is neglected. Approximately $10^6$ particles are used for the simulations, which is considered as medium to high resolution for free surface flows of these length scales. The simulations are performed on a single GPU-unit with 2048 cores and the average simulation time is approximately 14 hours. A typical time-dependent result for a height of 3m and a length of 7m, along with the perpendicular wall force component $F_x$, is shown in figure \ref{fig:Results3DExample}. The four consecutive instances of the example simulation show: initial liquid configuration; wave development due to pressure gradient; breaking wave impact on wall; liquid after impact. Depending on the height difference between both liquid columns, wave breaking can occur before the wave impacts the wall. In reality, a gas pocket is entrapped during this process, which is compressed and leads to a pressure build-up inside the gas pocket. However, the simulations performed here are free surface flows and the gas phase is not taken into account, so the physics in the entrapped gas pocket is ignored.
\begin{figure}[!h]
\centering
\includegraphics[width = \textwidth]{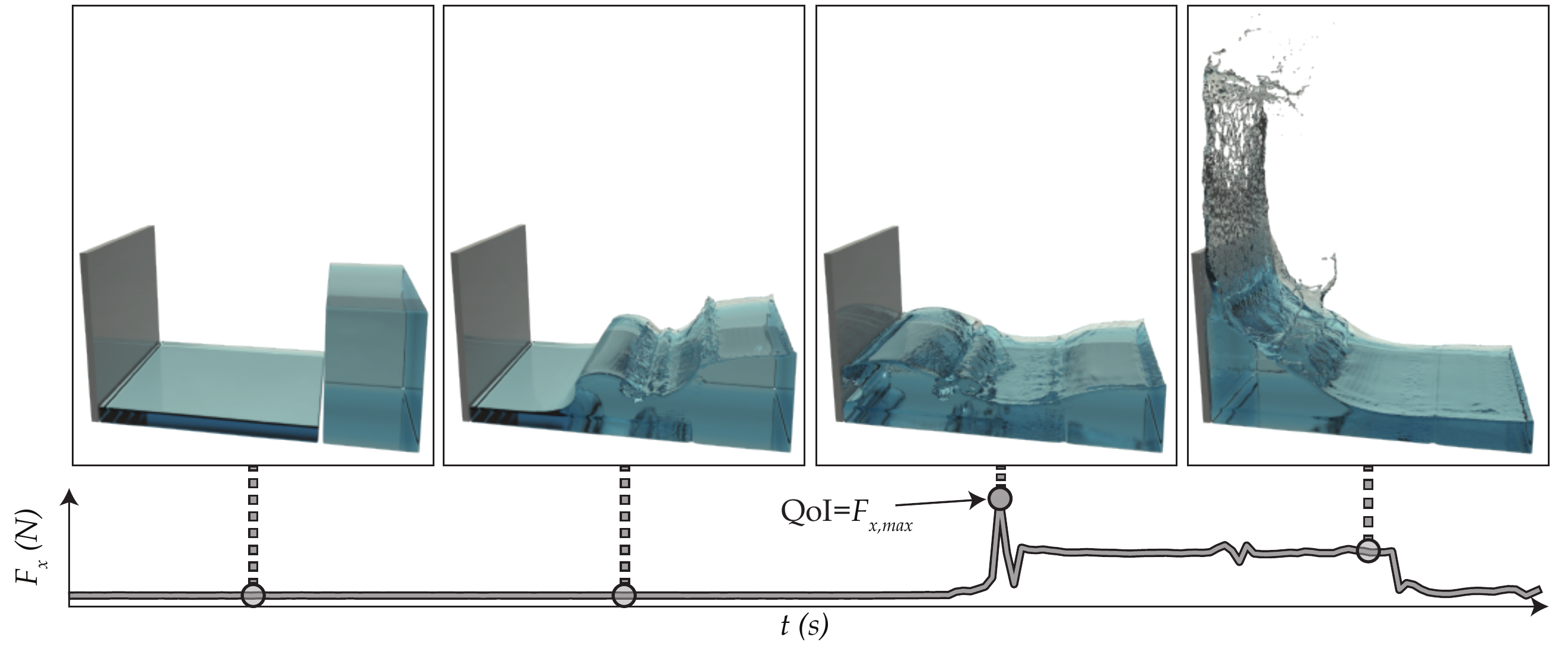}
\caption{\label{fig:Results3DExample} Single simulation with height of 3m and a length of 7m. The QoI is the maximum of the wall force $F_x$.}
\end{figure}\\
Since obtaining a parametric solution is the target, both parameters are assumed to be uniformly distributed. The MST-ME method is used to construct the unknown QoI response. A total of 6 iterations of the sampling algorithm is performed with weight function \eqref{eq:weight3}, which results in 30 samples. We expect the MST-ME method to automatically distinguish smooth and discontinuous behaviour of the QoI, which is important for this problem, since we have limited initial knowledge of the QoI as a function of the parameters. The results are shown in figure \ref{fig:Results3DResponse}.
\begin{figure}[!h]
\centering
\includegraphics[width = \textwidth]{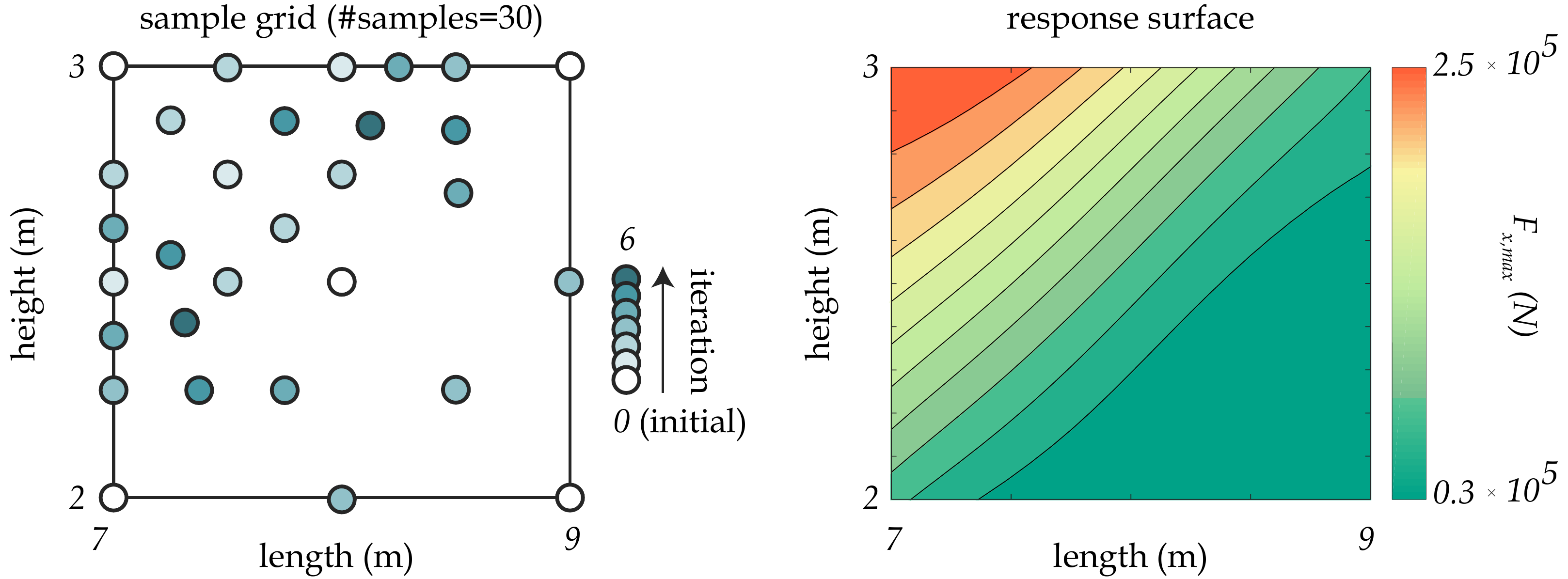}
\caption{\label{fig:Results3DResponse} QoI of the 3D dam break problem, obtained by the MST-ME method with 6 iterations.}
\end{figure}\\
The results indicate a smooth QoI response, which is approximately constant along the lines height/length=constant. This implies that for free surface flow, the force on the wall depends roughly on the ratio of height and length, and not on their separate values. The wall force increases when this ratio increases, which is intuitive from a physical point of view. Figure \ref{fig:Results3DResponse} (right) shows that this increase is non-linear, which corresponds to results previously reported for dry-bed dam break problems \cite{lobovsky_experimental_2014}.

Interestingly, in contrast to the SWE test-case, the QoI shows no discontinuity in the parameter space. This is possibly due to neglecting the gas phase in the free surface simulations. When simulated with a gas phase, the pressure build-up in the entrapped gas pocket (see figure \ref{fig:Results3DExample}) may lead to a discontinuity in the QoI. The strength of our proposed MST-ME method is that we do not require knowledge about the characteristics of the QoI beforehand, as it distinguishes smooth and discontinuous behaviour automatically.
\section{Conclusion}
\noindent In this paper we have presented a novel domain decomposition based interpolation method, the Minimum Spanning Tree Multi-Element (MST-ME) method. The unique property of the MST-ME method is that it adaptively constructs a surrogate model as a function of a set of uncertain
parameters for both smooth and discontinuous quantities of interest. The three key ingredients in the MST-ME method are: (I) adaptive sampling based on a minimum spanning tree with a smart weight function, (II) discontinuity detection and sample classification with a support vector machine algorithm, and (III) least orthogonal interpolation to construct local approximations. This combination of robust methods makes the MST-ME method a very practical method that is applicable to a wide range of UQ problems.

The MST-ME method has been applied to several numerical examples: domain decomposition, Genz function approximation, 1D shallow water equations, and a 3D dam break problem. Any discontinuities present in these test-cases are effectively captured by the method. In all cases, fast convergence is obtained, leading to an accurate surrogate model already at a tractable number of model runs. This surrogate model can be directly used as a fast tool for uncertainty quantification (for example with Monte Carlo type methods), but it is also a great tool for the parametric solution of black-box models, including partial differential equations. We also foresee application of this surrogate model in the solution of inverse problems. The freedom in the weighing function of the minimum spanning tree offers many applications, such as adaptive sampling for reliability analysis, where the weighing function may be adapted such that samples are placed in regions of low probability.

Currently, the MST-ME method does not include the option to detect discontinuities that run partially through the random space. Furthermore, the MST-ME method can not yet preserve the symmetry in the node distributions, which might be advantageous in certain special cases (e.g. when both the model and the underlying PDF of the random variables are symmetric). Lastly, the least orthogonal interpolation does not necessarily use all sample points in the construction of the local approximation. By adding a term to the weight function of the minimum spanning tree, which accounts for the stability of the interpolant (as is done for example in Leja nodes \cite{narayan_adaptive_2014}), the sample locations may be further improved, such that all samples are used in the construction of the interpolant.

\section*{Acknowledgements}
\noindent This work is part of the research programme ''SLING'' (Sloshing of Liquefied Natural Gas), which is (partly) financed by the Netherlands Organisation for Scientific Research (NWO).
\section*{References}
\bibliography{MyLibrary}
\end{document}